\documentclass[10pt,twoside]{article}
\usepackage{amsmath}
\usepackage{epsfig, color, colordvi}
\usepackage{verbatim,amssymb,amsmath}
\usepackage{tikz}

\newcommand{\qzero}{\hat q^2_0}
\newcommand{\cstar}{c^2_*}
\newcommand{\phisup}{\varphi^{sup}}
\newcommand{\psup}{p^{sup}}

\newtheorem{theorem}{Theorem}[section]

\newtheorem{remark}[theorem]{Remark}
\def\QED{~\rule[-1pt]{8pt}{8pt}\par\medskip}

\begin{document}
\title{ A  variational inequality formulation for transonic compressible steady potential flows: \\
Radially symmetric transonic shock
}
\author{Yung-Sze Choi\thanks{
Department of Mathematics,
University of Connecticut,
Storrs, CT 06269-3009, 
email: choi@math.uconn.edu.
}
\and
Eun Heui Kim
\thanks{
Department of Mathematics,
California State University,
Long Beach, CA 90840-1001,
email: EunHeui.Kim@csulb.edu.
Research supported by the National Science Foundation under the grant DMS-1615266.
}
}
\date{}
\maketitle

\begin{abstract}
 We establish {a} variational inequality formulation that captures the transonic shock for {a} steady compressible potential flow.
{Its critical point satisfies the transonic equation; moreover the associated jump conditions across its free boundary
 match the usual Rankine-Hugoniot jump conditions for a shock.}   By means of example we validate our formulation,
{and establish} the necessary and sufficient condition for the existence of a transonic shock. Numerical results are also discussed.

\noindent
AMS: {Primary: 76L05, 35L65; Secondary: 65M06, 35M33.}
\end{abstract}


\section{Introduction}\label{intro}
\setcounter{equation}{0}
\renewcommand{\theequation}{\thesection.\arabic{equation}}

In an effort to understand steady transonic {potential} flow in compressible gas dynamics we introduce {a} variational inequality formulation.
To validate our formulation, we build a model problem with an exact solution, 
and {establish} 
the necessary and sufficient conditions for the existence of transonic shock in our model using the new variational scheme.
While we study a simple model system at this point, 
{our variational inequality formulation is applicable to more general multidimensional settings.}

{It is well known that the steady compressible Euler equation is the governing model system for various applications such as a flow over a profile (a steadily moving aircraft). Consequently understanding boundary value problem in multidimensional gas dynamics has always been of interest to mathematicians and engineers alike.}
The classical book by Courant and Friedrichs \cite{CoFr} contains many more such applications that involve
steady compressible Euler equations in multidimensional settings.
The following quoted from Morawetz \cite{Morawetz} has been a long standing open problem:
\begin{quote} 
In the study of transonic flow, one of the most illuminating theorems to prove would be:

{\em Given an airfoil profile and a continuous two dimensional irrotational transonic compressible inviscid flow past it with some given speed at infinity, there does not exist a corresponding flow with a slightly different speed at infinity.}
\end{quote}
In order to answer Morawetz's question, one must face issues arising from the multi-dimensionality of the boundary value problem.
For a given subsonic constant speed at infinity, Shiffman \cite{Shiffman} used the direct method of calculus of variations 
to establish an existence theorem for the subsonic flow in the entire region around a profile. 
The existence and uniqueness results for this subsonic flow over the profile is completed by Bers \cite{Bers}, a stronger uniqueness result by Finn and Gilbarg \cite{FG} and for higher dimensions by Dong and Ou \cite{DongOu}. 
All these results rely upon the condition that the {given} speed at infinity must be {\em sufficiently} subsonic.   
{To} our best knowledge, 
it is unknown whether a critical point for the variational problem exists for any arbitrary speed at infinity. 
 
In fact subsonic flow {at} infinity does not necessarily remain subsonic everywhere; it can {become} transonic and create a shock over a convex profile \cite{CoFr}.
For simpler configurations when the supersonic flow at infinity is given, a shock wave is formed in front of the profile. 
In general the shock location is not known apriori. 
{In other words a steady transonic flow that contains a shock is a free boundary problem.}

Recent progress on transonic problems reports on specific configurations (many are on regular shock reflections with the known supersonic state){,} specific model systems including
 {potential} flow \cite{ChFe,Chen,EL,potell, Mor, XY}, the steady transonic small disturbance (STD) equations \cite{CKL},  the unsteady transonic small disturbance (UTSD) equations \cite{bh1,CK1,CKK:WRR,CKL},
the nonlinear wave equations \cite{CKK:nlwe, JKC, sever},
the pressure gradient equations \cite{KS,Zheng1}, and references therein.
{These results rely} upon perturbation methods,  partial hodograph methods, and various fixed point iteration methods.

 We propose to tackle the transonic problem using a variational inequality 
 formulation with  the  shock location represented by
 a free boundary to be determined as a part of the solution. 
 Br\'ezis and Stampacchina \cite{BreStam} showed that the problem of infinite plane flow around a symmetric convex profile can be reduced to a variational inequality in the hodograph plane, and established the existence of a subsonic flow.
  Shimborsky \cite{Shi} utilized {the} variational inequality in the hodograph plane setting studied in \cite{BreStam} to establish a subsonic symmetric channel flow (Laval nozzles).
  Both \cite{BreStam,Shi} results focused on subsonic flow.
Later a few attempts have been tried to extend the
 formulation to cover supersonic flow as well as shock formation
\cite{FMN, Gittel} and the references therein. For example
 in \cite{Gittel}, an artificial entropy condition was imposed on admissible
 functions so as to ensure the existence of a minimizer to a functional. 
 We note, however, that it is not clear a priori whether the minimizer is a solution of the transonic problem.
 Furthermore, their variational formulations do not provide any information how to locate a shock,
 and whether the Rankine-Hugoniot jump conditions are satisfied. 
 
 In this paper, we consider a configuration for the 
given supersonic upstream flow be isentropic. Since the disturbance
caused by {an immersed} body can only propagate downstream, the supersonic flow will
not be affected {by the formation of a} shock. 
In general the flow will not be isentropic behind the shock
(it can be, if the entropy  increases by the same amount across the shock).
 In other words
 the replacement of the full system governing compressible fluid flow
by the transonic equation induces inconsistency. Such an inconsistency
 usually {results} in not satisfying all the three Rankine-Hugoniot conditions
 across a shock.
 
The main goals of this paper is to establish a correct variational inequality formulation 
incorporating a transonic free boundary,
and find the corresponding {critical point} 
which satisfies the Rankine-Hugoniot jump conditions across the transonic shock.
More precisely, 
we consider a configuration for which the supersonic flow upstream of the shock is known,
so that one can find a potential function, {denoted by} $\varphi^{sup}$, which gives rise to this
simple flow. We next consider
 the shock as a free boundary where the supersonic flow in the direction
 normal to the shock makes an abrupt
 change to a  subsonic flow. Since the flow behind the shock is assumed to be potential, 
 it suffices to find the corresponding subsonic potential function $\varphi^{sub}$, 
 satisfying $\varphi^{sub} \geq \varphi^{sup}$ with $\varphi^{sub}=\varphi^{sup}$
 at the shock.  This obstacle problem is a free boundary
 problem that can be formulated in a variational inequality setting. 
 That is, we perform the variational calculus on a closed convex set in
  a function space
rather than on the whole function space. 
 Such a formulation 
will automatically guarantee the continuity of tangential velocity across a shock,
which is a physical requirement (see equation (81.6), p.318, \cite{Landau}).
The remaining hurdle is to ensure the natural boundary condition associated with
the variational functional agrees with some (or a linear combination) of the
Rankine-Hugoniot conditions. 
As remarked earlier one cannot expect all such conditions to be satisfied as the
 flow behind the shock may not be isentropic. 
 
To test the correctness of our variational inequality formulation
we consider a model problem that involves only radial flow. 
This problem admits an exact solution which can be used to validate our formulation. With 
  the actual flow remains potential behind the shock, it is remarkable that our 
  variational formulation
 captures all the Rankine-Hugoniot conditions in this case 
 so that its critical point, which solves the transonic equation, is also
 an exact solution to the full system of compressible fluid flow with shock wave.
 This is only possible since in our formulation 
 we have allowed entropy to change behind the shock.
 We note, however, that
 there are many works that assume the same constant entropy 
 in front and behind the shock for simplicity.
 This automatically eliminates any chance of capturing the exact solution to
 the full system of compressible fluid flow in variational formulation.

 More importantly, the model problem shows
 clearly that the critical point of our variational inequality formulation corresponds to 
a shock solution to the transonic equation is a {\em saddle} point. 
Our variational formulation ensures that the saddle point
satisfies one Rankine-Hugoniot condition and a linear combination of the remaining two. Thus one can think of it satisfying two
such conditions. And if the flow is actually potential behind the shock,
like in the model problem, one can have all three jump conditions being satisfied.

   The paper is organized as follows.
  In section~\ref{sec_gov} we present an overview of the steady flow of compressible system. 
  We pay a particular attention to the physical entropy constant which is typically normalized and simplified elsewhere.
  Section~\ref{sec_model} comprises our model problem and its exact solution, see Subsection~\ref{sec_exact}. This is also substantiated by 
 numerical results. 
 In Section~\ref{sec_variation}, we discuss the variational inequality formulation and the free (transonic shock) boundary condition. 
 We show that {its critical point} satisfies the mass conservation equation and {some} jump conditions on the free boundary that 
 {corresponds} to the transonic shock.
 We establish the existence of the saddle point; {the necessary and sufficient conditions of such an existence}  agree with that obtained in Section~\ref{sec_exact}.
 We conclude the paper by presenting numerical results that demonstrate the transonic solution being the saddle node of the variational function in Section~\ref{sec_num}.
 
 We believe our results will serve as a vehicle for understanding transonic flows in particular the long standing open problem of the flow over a profile.

\section{Review of governing equations} \label{sec_gov}
\setcounter{equation}{0}

{Conservation laws} are governing principles in gas dynamics.
While there are many references dedicated to explain compressible flow, 
we give a short survey on the conservation principles which will be useful for our variational formulation.

The two dimensional steady {state} compressible Euler system for an ideal gas in conservation form reads
\begin{eqnarray}\label{massrec}
(\rho u_1)_x + (\rho u_2)_y &=&0, \quad {Conservation \ of \ Mass}\\
\label{mom1}
(\rho u_1^2 +p)_x + (\rho u_1 u_2)_y&=&0, \quad {Conservation \ of \ Momentum}\\
(\rho u_1 u_2)_x + (\rho u_2^2 +p)_y &=&0,\label{mom2} \\\label{enrrec}
(\rho u_1(\frac{1}{2} q^2 +\iota))_x + (\rho u_2 (\frac{1}{2} q^2 + \iota))_y&=&0.\quad {Conservation \ of \ Energy}
\end{eqnarray}
Here $\rho$ is the density, $p$ is the pressure,  ${\bf u}=(u_1,\ u_2)$ is the velocity 
vector, $q=\sqrt{u_1^2+u_2^2}$ is the gas speed and $\iota$ is the enthalpy per unit mass 
of the gas. 

Entropy of any ideal gas particle remains constant
during its motion except when crossing a shock.
As a consequence one has
\begin{equation} \label{polytropic}
p=k \rho^{\gamma}
\end{equation} 
 for some positive constants $1<\gamma<\infty$ (typically $1< \gamma \le 5/3$) 
 and $k$.
Physically, the adiabatic exponent $\gamma$ is the ratio of specific heat capacity per unit mass
 at constant pressure to that at constant temperature for the gas, 
and $k$
is a {parameter} that depends {only} on the local entropy \cite[p.6]{CoFr}. In case when
{upstream flow is homogeneous,
the entropy, and hence $k$, are
the same everywhere in this spatial region.}
{After} a gas particle crosses a shock, {its entropy increases and so does the parameter $k$.} 
{With (\ref{polytropic})} the enthalpy of an ideal gas reads
 \begin{equation} \label{def_i}
\iota= \frac{\gamma}{\gamma-1} k \rho^{\gamma-1} = \frac{\gamma}{\gamma-1} \frac{p}{\rho} =\frac{c^2}{\gamma-1},
 \end{equation}
 where $c$ is the (local) sound speed defined by the above equation.

We consider a potential flow that is supersonic far upstream; {the potential function is denoted by $\varphi^{sup}$ so that}
${\bf u}=- \nabla \varphi^{sup}$. 
The negative sign in the last equation {is convenient} in our variational inequality formulation that will be discussed later. 

The {subsequent} steady adiabatic flow satisfies the Bernoulli's law
\begin{eqnarray}\label{B-law}
\frac{1}{2}q^2 + \frac{c^2}{\gamma-1} = \frac{1}{2} \hat q^2_0 \;,
\end{eqnarray}
 where $q$ is the speed of the flow and $\hat q_0$ is a {known positive}
 constant along a streamline determined by the given upstream flow. {In general}
 the Bernoulli's constant $\frac{1}{2} \hat q^2_0$ may have different values along different streamlines, and the same is true for its entropy.
 {Since the upstream flow is homogeneous in our case,}
$\frac{1}{2} \hat q^2_0$ is the same on every streamline. 

{Where} a stationary shock is formed,
one obtains from (\ref{massrec})-(\ref{enrrec})
  the Rankine-Hugoniot jump conditions 
\begin{eqnarray}\label{sk1}
[\rho v]&=&0\;, \\\label{sk2}
[\rho v^2 +p]&=&0\;, \\\label{sk3}
[\rho v(\frac{1}{2} v^2 + i)]&=&0\;.
\end{eqnarray}
Here $[ \cdot ]$ in the above equations denotes the jump across the shock
and $v$ is the velocity component normal to the shock. Thus there is 
no jump in $v^2/2+\iota$ across a shock.  Since the tangential speed across
a shock is continuous (which can be considered as a $4^{th}$ Rankine-Hugoniot condition),
this implies $q^2/2+\iota\;$ is continuous across a shock and hence
equals to the same constant $\hat q_0^2/2$ everywhere immediate behind a
shock. 

After crossing the shock, the flow may no longer be potential. {As}
the Bernoulli's equation requires $q^2/2+\iota$
be constant along each individual stream line for a steady state,  one can conclude
that (\ref{B-law}) holds everywhere before and after the shock {with the same}
Bernoulli's constant $\hat q_0^2/2$ {as in upstream flow}. As a consequence, we can replace
 the energy equation (\ref{enrrec}) and its jump condition (\ref{sk3}) by (\ref{B-law}).

We now adopt the usual simplification  that the flow behind a shock is 
potential with its velocity vector
${\bf u}^{sub}=- \nabla \varphi^{sub}$ for some
velocity potential function $\varphi^{sub}$ to be determined. Even the shock is a strong
one, this may still be a good approximation so long as the increase in entropy
across the shock is more or less uniform, as illustrated by the model problem below.
We like to emphasize that even we use the superscript $\;sub\;$ to designate the 
variables behind a shock, the flow may not be necessarily subsonic everywhere.

From the Bernoulli's law (\ref{B-law}) and the flow being assumed to be potential, we can now {\em{define}}
both $c^2$ and $\rho$ as functions of the speed $|\nabla \varphi|$ everywhere
\begin{eqnarray}\label{soundspeed}
c^2   \equiv \frac{\gamma-1}{2} ( \qzero - |\nabla \varphi|^2)= k \gamma \rho^{\gamma-1}, 
\end{eqnarray}
and
\begin{eqnarray}\label{rhov}
\rho =\rho(|\nabla \varphi|^2) &:= & \left(\frac{1}{k} \frac{\gamma-1}{2\gamma}(\qzero -|\nabla\varphi|^2)\right)^{\frac{1}{\gamma-1}} \;
\end{eqnarray}
with $\varphi$ being $\varphi^{sup}$ or $\varphi^{sub}$. 
The entropy increase across the shock results the changes in  
$k$ that increases from $k^{sup}$ to a larger value $k^{sub}$ after crossing
the shock.  It is important to note a major difference between
(\ref{B-law}) and (\ref{rhov}): 
the potential flow assumption has been built in the latter equation already.
{If the flow behind the shock were truly potential, the momentum equations will be automatically
satisfied.}

Consequently,
the continuity equation (\ref{massrec}) behind the shock
can be written in an equation of a potential flow $\varphi$
\begin{eqnarray}\label{mass_div}
div( \rho\ \nabla \varphi) &=& div( \rho(|\nabla \varphi|^2)\ \nabla \varphi) =0.
\end{eqnarray}
where $\rho$ satisfies \eqref{rhov}.
This continuity equation can be written in a non-divergence form, and denoted by,
\begin{eqnarray}\label{pde}
Q\varphi := (c^2 -\varphi_x^2 ) \varphi_{xx} -2\varphi_x\varphi_y \varphi_{xy} +(c^2 -\varphi_y^2) \varphi_{yy} =0.
\end{eqnarray}
In fact $Q\varphi = \sum_{i,j} A_{ij} \, D_{ij} \varphi$ with 
$A=\left( \begin{array}{cc} c^2-\varphi_x^2 & -\varphi_x \varphi_y \\
 -\varphi_x \varphi_y &  c^2-\varphi_y^2 \end{array} \right)$. It is readily verified that
\begin{eqnarray}\label{det}
det(A) &=& c^2 (c^2 - \varphi_x^2 -\varphi_y^2)
\end{eqnarray}
and thus the operator $Q$ is elliptic whenever $c^2 > |\nabla \varphi|^2$, that is, 
\[ \cstar \equiv \frac{\gamma-1}{\gamma+1}\qzero > |\nabla\varphi|^2 \]
which is equivalent to subsonic flow.

Our model problem, {which will be used to validate our variational inequality formulation},
 concerns a radial fluid flow without angular motion. 
Let $r$ be the radial distance from the origin and
and $v$ denote the radial velocity (positive in the increasing $r$ direction)
which is a function of $r$ only.
It is convenient to write the Euler system in such a coordinate system:
\begin{eqnarray}\label{mass}
(r\rho v )_r &=&0, \quad {Conservation \ of \ Mass}\\\label{mom}
(r\rho v^2 +rp)_r &=&p,\quad {Conservation \ of \ Momentum}\\\label{enr}
(r\rho v (\frac{1}{2} v^2 + \iota))_r &=&0.\quad {Conservation \ of \ Energy}
\end{eqnarray}
With the assumption of potential flow behind the shock, one {only} needs to study
(\ref{mass}) and  (\ref{rhov}). 
Since $\frac{d\iota}{dr}=\frac{1}{\rho} \frac{dp}{dr}$,
it is easy to verify the momentum equation (\ref{mom}) is an easy consequence of 
(\ref{rhov}). Hence the solution to the transonic equation in this case
turns out to be exact for the Euler system, because radial flow behind a shock is always potential.

\section{A model problem}  \label{sec_model}
\setcounter{equation}{0}

Let $a, R$ be given positive constants with $a<R$. 
We consider a simple configuration that involves radial flow in the annulus
$\{ \mathbf{x} \in \mathbf{R}^2: a<|\mathbf{x}|<R \}$ with 
 a circular shock inside.
We study \eqref{mass} and \eqref{rhov} under the assumption that the flow is potential behind the shock.
 Since radial flow is always potential, a solution to these two equations
 will be the exact solution for the Euler system (\ref{mass})-(\ref{enr}).
Impose the boundary conditions
\begin{eqnarray}\label{bcR}
v=-v_0,\; \rho=\rho_0,\;  k=k_0 & {\rm at } & r=R \;,\\\label{bca}
v=-v_a  & {\rm at } & r=a\;,
\end{eqnarray}
with both $v_0, v_a>0$.
The given velocities  at both $r=a$ and $r=R$ are negative since they point towards the origin.
Because of (\ref{mass}) and (\ref{sk1}), the total flux is a constant everywhere,
regardless to be in front or behind the shock.
Combining this information with (\ref{rhov}) and
 the boundary conditions at $r=R$, we obtain
\begin{eqnarray*}
- \rho v r &=& \rho_0 v_0 R := M_0 \;, \\
 \frac{1}{2} v^2 + \frac{\gamma}{\gamma-1} k\rho^{\gamma-1} &=&  
 \frac{1}{2} v_0^2 + \frac{\gamma}{\gamma-1} k_0 \rho_0^{\gamma-1} := \frac{\qzero}{2}.
\end{eqnarray*}
That is, $v$ and $\rho$ can be written in terms of $r$ and $k$ by using
\begin{eqnarray}\label{soln}
- \rho v = \frac{M_0}{r},\quad 
\rho= \left(\frac{1}{k} \frac{\gamma-1}{2\gamma}(\qzero -v^2)\right)^{\frac{1}{\gamma-1}}
\;.
\end{eqnarray}
We let the boundary conditions (\ref{bcR}) and (\ref{bca}) satisfy
\begin{eqnarray} \label{supR}
c^2|_{r=R}= \gamma k_0\rho_0^{\gamma-1} &<& v_0^2, \quad { Supersonic\ at\ r=R }\;,\\
\cstar := \frac{\gamma-1}{\gamma+1} \qzero & >& v_a^2, \quad { Subsonic\ at\ r=a }\;,  \label{suba}
\end{eqnarray}
so that
\begin{equation}
v_a <c_* < v_0 < \hat{q}_0 \;.
\end{equation}
Employing (\ref{soln}) and the boundary condition (\ref{bca}), we conclude that
\begin{eqnarray}
k_a&=& 
\frac{\gamma -1}{2\gamma} \left(\frac{a v_a}{M_0}\right)^{\gamma-1}(\qzero - v_a^2),
\label{k_a}
\\
\rho_a & = &\frac{M_0}{a v_a}. \label{rho_a}
\end{eqnarray}
We note that $k_a$ depends only on the given boundary conditions, but not
on any other features of the subsonic solution behind the shock. 
 One expects the change in entropy $k$ may, in general, depends on the location of the shock.

\subsection{Exact solutions to the model system} \label{sec_exact}

The radially symmetric model has a stationary shock which satisfies all the Rankine-Hugoniot conditions
 under appropriate conditions, and the exact solution {will} be found explicitly; {in particular we need
 to pin down the shock location.}
  Later in Section~\ref{sec_analysis1D},  we obtain a solution from our variational formulation, which is in exact match with the explicit solution obtained in this section. This validates our methodology. 
 
 For {weak shock} it is known that the change in entropy is of third order {of the pressure change} \cite[p.323]{Landau}.
 {This is the rationale behind many transonic potential flow studies for setting the entropy constant $k$ to stay the same after passing the shock.}
However {we test our
our variational inequality formulation without the weak shock assumption in the model.}
To follow the change in $k$,
from (\ref{soln}) we have
 \begin{equation} \label{eqnK}
 k= 
\frac{\gamma -1}{2\gamma} \left(\frac{-r v}{M_0}\right)^{\gamma-1}(\qzero - v^2).
\end{equation}

\begin{figure}[ht]
\begin{center}
\includegraphics[height = 1.8in,width = 2.3in]{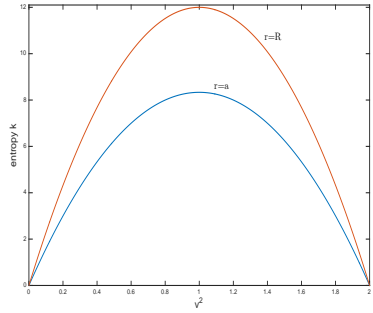}
\end{center}
\vspace*{-4mm}
\caption{$k$ curves: Matlab plot with $\gamma=3$, $\cstar=1$, $M_0=1$, $R=6$ and $a=5$.}
\label{kcurves}
\end{figure}
\noindent
Figure~\ref{kcurves} depicts the $k$ plots with respect to the speed $v^2$ where $r$ is at $r=R$ and $r=a$.  
 A maximum point of the graph $k$ is always located at $|v|=c_*$; {in this example $c_*=1$}. 
 When $r=a$,  the corresponding
 \begin{equation} \label{kmax}
k_{max,\, r=a}=\frac{1}{\gamma} ( \frac{a}{M_0})^{\gamma-1}  c_*^{\gamma+1}\;.
\end{equation}
To satisfy
the increase of entropy after the shock, we need $k_a>k_0$. Therefore 
a necessary condition for a shock wave to exist is
\begin{equation} \label{k_condition}
k_0< k_a< \frac{1}{\gamma} ( \frac{a}{M_0})^{\gamma-1}  c_*^{\gamma+1}\;.
\end{equation}
While the second inequality {is immediate}, 
in order to have a transonic flow, 
we choose a subsonic $v_a$ to be large enough so that $k_a>k_0$.

In the supersonic region in front of the shock, 
we have $k=k_0$. With any given $r>a$, one can solve for a unique supersonic $v$ from 
(\ref{eqnK}) and write
$-v= V^{sup} (k_0,r)$ (recall that $v<0$ so that $V^{sup}$ is a positive
function). Correspondingly we find $\rho=\rho^{sup}(k_0,r)$ by substituting $k=k_0$
and $v=-V^{sup}(k_0,r)$ in the second equation in (\ref{soln}).
Similarly in the subsonic region behind the shock, we have $k=k_a$ and
we can solve for a unique subsonic $-v$ so that
$-v=V^{sub}(k_a,r)$, and
 $\rho=\rho^{sub}(k_a,r)$ for any $r>a$.
 
Across the transonic shock (if it exists), denoted by $\Sigma$,
the Rankine-Hugoniot conditions take the form of
\begin{eqnarray}\label{rh1}
\rho_1 v_1 &=& \rho_2 v_2\\\label{rh2}
\rho_1 v_1^2 + k_0\rho_1^\gamma&=& \rho_2 v_2^2 + k_a \rho_2^\gamma\\\label{rh3}
 \frac{1}{2} v_1^2 + \frac{\gamma}{\gamma-1} k_0\rho_1^{\gamma-1}&=& \frac{1}{2} v_2^2 + \frac{\gamma}{\gamma-1} k_a\rho_2^{\gamma-1}
\end{eqnarray}
where $\rho_1=\rho^{sup}(k_0, r_s), - v_1=V^{sup}(k_0, r_s)$ and $\rho_2=\rho
^{sub}(k_a, r_s), -v_2=V^{sub}(k_a, r_s)$.
Note that (\ref{rh1}) and (\ref{rh3}) are satisfied because of the first equation in (\ref{soln}) and 
our choice of Bernoulli's constant being unchanged.
If there exists $r=r_s$ satisfying \eqref{rh2}, we have a transonic shock conserving mass, momentum and energy across it.
We find the transonic shock that conserves the momentum by including the changes in entropy, that is the entropy constant $k$ is changed across the shock.\footnote{For the shock conditions for potential flow, Morawetz \cite{Mor} noted that the momentum, just as the entropy and vorticity, change is of third order in the shock strength, and stated further that one must give up conservation of momentum normal to the shock.}

We observe that the momentum flux can be written as
\begin{eqnarray*}
\rho v^2 + p &=& \rho (v^2 + \frac{1}{\gamma} c^2)\\
&=& \rho (v^2 +\frac{\gamma-1}{2\gamma} (\qzero - v^2))\\
&=& \rho (\frac{\gamma-1}{2\gamma} \qzero + \frac{\gamma+1}{2\gamma} v^2)\\
&=&\frac{\gamma+1}{2\gamma} \rho (\cstar + v^2)\\
&=& \frac{\gamma+1}{2\gamma} \frac{M_0}{r} (\frac{\cstar}{-v} + (-v)).
\end{eqnarray*}
Let us denote
\begin{eqnarray} \label{def_H}
H&:= & \frac{\gamma+1}{2\gamma} M_0 (\frac{\cstar}{-v} + (-v)) = r(\rho v^2 +p).
\end{eqnarray}
Recall that $-v=V^{sup}(k_0,r)$ in front of the shock and $-v=V^{sub}(k_a,r)$ behind the shock.
We then consider $H$ as a function of $r$:
\[
H^{sup/sub}(r) = \frac{M_0(\gamma+1)}{2\gamma} (\frac{\cstar}{V^{sup/sub}} + V^{sup/sub})\;,
\]
where we  have,
and sometimes will again, suppressed the dependency of $H$ on $k$ for notational simplicity.
In the $(r, H)$ coordinates we let point $A$ be $(a, H^{sub}(k_a,a))$ and 
point $C$ be $(R, H^{sup}(k_0,R))$, which can be computed easily
using the given boundary
conditions at $r=a$ and $r=R$, respectively.
\begin{figure}
\begin{center}
\includegraphics[height = 1.8in,width = 2.3in]{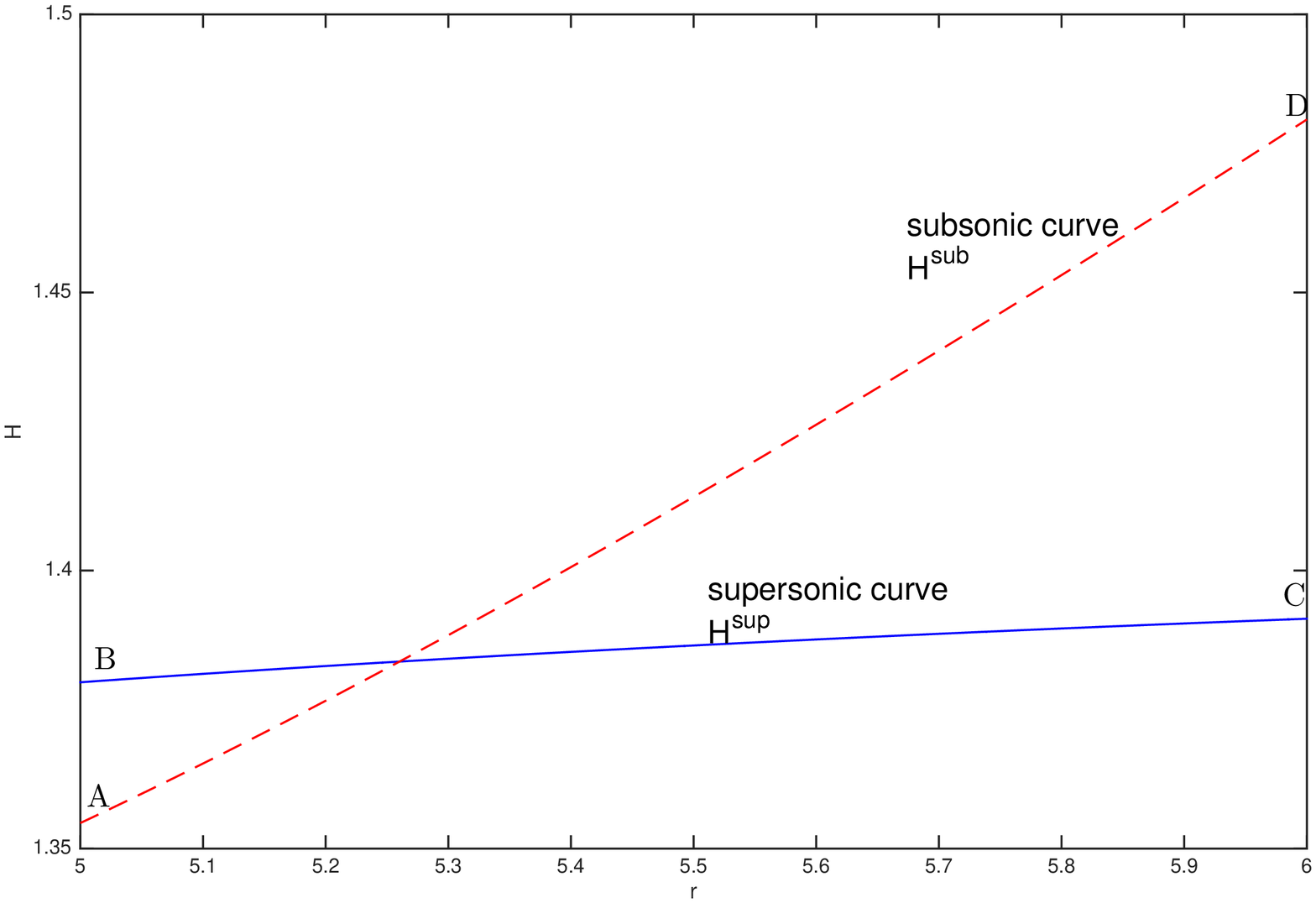}
\end{center}
\vspace*{-4mm}
\caption{$H$-curves: Matlab plot with $\gamma=3$, $\cstar=1$, $M_0=1$, $R=6$, $a=5$, $v_a= \sqrt{0.7} > \nu =\sqrt{0.5}$ and $v_0=\sqrt{1.8}$. The corresponding $k$ values are $k_a=7.58\overline{3}$ and $k_0=4.32$. $r_s=5.260220746$}
\label{Hcurves}
\end{figure}

Figure~\ref{Hcurves} depicts the curves $H^{sup}(k_0,r)$ and $H^{sub}(k_a,r)$ as functions of $r$.
Then there must be the corresponding end points $B=(a,H^{sup}(k_0,a))$ and $D=(R, H^{sub}(k_a,R))$, respectively.
Thus the shock position $r=r_s$ corresponds to a point at which the curves $H^{sup}$ and $H^{sub}$ cross. 
If they do cross, the intersection point has to be unique. The uniqueness 
can be easily seen as follows. From (\ref{soln}) and
$p= k \rho^{\gamma}$, one easily verifies that
\begin{equation} \label{pressure}
p=\frac{\gamma-1}{2 \gamma} \frac{M_0}{r} \frac{{\hat q_0^2}-|v|^2}
{|v|} 
\end{equation}
which is a decreasing function for $0 \leq |v| \leq \hat{q}_0^2$.
Since $V^{sup}>c_*> V^{sub}$, we have $p^{sub}>p^{sup}$. Thus
as a consequence of (\ref{mom}), we have $dH^{sub}/dr=p^{sub}>p^{sup}=dH^{sup}/dr$ 
for all $r$, which leads to a unique crossing point.

We now turn to the existence of the intersection of the curves. 
Let conditions  \eqref{supR}, \eqref{suba} and \eqref{k_condition} hold.
Since the curve $H^{sub}$ always has a steeper slope than the curve $H^{sup}$,
for given $A$ and $C$ the necessary and sufficient conditions for crossing are:
\begin{eqnarray} \label{cross_cond}
H_D= H^{sub}(k_a,R) & > & H^{sup}(k_0,R) =H_C \;, \\
H_B=H^{sup}(k_0,a) &> & H^{sub}(k_a,a)=H_A \;. 
\end{eqnarray}
These are equivalent to
\begin{eqnarray*}
 H_D = H^{sub}(R)& =&\frac{M_0(\gamma+1)}{2\gamma} (\frac{\cstar}{V^{sub}(R)} + V^{sub}(R)) \\
 &>&  \frac{M_0(\gamma+1)}{2\gamma} (\frac{\cstar}{V^{sup}(R)} + V^{sup}(R)) = H^{sup}(R) =H_C,
 \end{eqnarray*}
 and
\begin{eqnarray*}
 H_B = H^{sup}(a)& =&\frac{M_0(\gamma+1)}{2\gamma} (\frac{\cstar}{V^{sup}(a)} + V^{sup}(a)) \\
 &>&  \frac{M_0(\gamma+1)}{2\gamma} (\frac{\cstar}{V^{sub}(a)} + V^{sub}(a)) = H^{sub}(a) =H_A.
 \end{eqnarray*}
Further simplifications lead to
\begin{eqnarray}\label{cond1}
\cstar &>& V^{sub}(k_a,R)V^{sup}(k_0,R),\\
\cstar &<& V^{sup}(k_0,a) V^{sub}(k_a,a).\label{cond2}
\end{eqnarray}
Thus \eqref{cond1} and \eqref{cond2} are the necessary and sufficient
conditions for the curves $H^{sup}(k_0,r)$ and $H^{sub}(k_a,r)$ to cross
at $r=r_s$ with $a<r_s<R$. Once we know the position of the shock, the
exact solution is completely known and satisfies all the Rankine-Hugoniot
conditions and the Euler system \eqref{mass}-\eqref{enr}.

We point out that conditions \eqref{supR}, \eqref{suba}, \eqref{k_condition} 
do not imply \eqref{cond1} and \eqref{cond2}.
In fact using the Rankine-Hugoniot jump conditions,  it is shown in \cite[(67.03), p.148]{CoFr} that the density compression ratio $\rho_2/\rho_1$
due to a shock in an ideal gas is always restricted to
\begin{eqnarray*}
\frac{\gamma-1}{\gamma+1} < \frac{\rho_2}{\rho_1}& = &\frac{V_1}{V_2}< 
\frac{\gamma+1}{\gamma-1},
\end{eqnarray*}
where we denote $V_1=V^{sup}(k_0, r_s)$ and $V_2=V^{sub}(k_a, r_s)$,
and consequently
\begin{eqnarray}\label{shkcomp}
\frac{\gamma-1}{\gamma+1} V_1&<& V_2 \;. 
\end{eqnarray}
Employing the Prandtl's relation \cite[page 147]{CoFr} at the shock, i.e.
\begin{eqnarray*}
\cstar &=& V_1 V_2,
\end{eqnarray*}
the above inequality \eqref{shkcomp} becomes
\begin{eqnarray*}
v_2^2 = V_2^2 &> & \cstar \frac{\gamma-1}{\gamma+1}, 
\end{eqnarray*}
which implies yet another constraint on $v_a$,
\begin{eqnarray}\label{vacond}
v_a \ge V_2 &>& c_*\sqrt{\frac{\gamma-1}{\gamma+1}} := \nu.
\end{eqnarray}
Here we have used $|v|$ being a decreasing function of $r$ in the subsonic
regime behind the shock. This can be easily seen from Figure~\ref{kcurves}
by observing the horizontal line $k=k_a$ intersects the curves at successively
smaller subsonic $|v|$ with increasing $r$. There exists $v_a$ satisfying
\eqref{supR}, \eqref{suba} and \eqref{k_condition}, but not \eqref{vacond}.
We emphasize that condition \eqref{vacond} is still a necessary condition
for the $H^{sup}$ and $H^{sub}$ curves to cross. The necessary and sufficient
conditions are (\ref{cond1}) and (\ref{cond2}).

\begin{figure}
\begin{center}
\includegraphics[height = 1.8in,width = 2.3in]{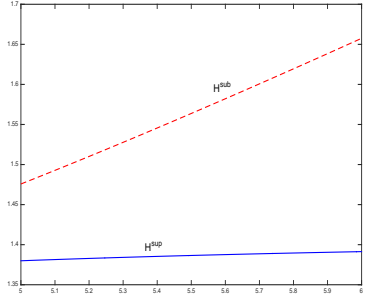}
\end{center}
\vspace*{-4mm}
\caption{A counter example of $H$-curves: Matlab plot with $\gamma=3$, $\cstar=1$, $M_0=1$, $R=6$, $a=5$, $v_a= \sqrt{0.4} < \nu =\sqrt{0.5}$ and $v_0=\sqrt{1.8}$. The corresponding $k$ values are $k_a = 5.\overline{3}$ and $k_0=4.32$.}
\label{Hcountereg}
\end{figure}

While data leading to Figure~\ref{Hcurves} shows that the H-curves cross,
Figure~\ref{Hcountereg} depicts a counter example in which $v_a$ satisfies all
the constraints except for \eqref{vacond}. More precisely we have
 $v_a =\sqrt{0.4} < \nu =\sqrt{0.5}$.
 As predicted the H-curves do not cross.
 
\section{Variational inequality formulation} \label{sec_variation}
\setcounter{equation}{0}
We now introduce a variational inequality formulation to solve the one-dimensional model problem
and show that the solution of the variational inequality is the same as the exact solution we found in Subsection~\ref{sec_exact}. The formulation is first derived for a general multidimensional case, and then reduced to the model problem of the radial flow.

Since the upstream supersonic flow is potential, there exists a $\varphi^{sup}$ such
that the velocity  $\mathbf{u}=- \nabla \varphi^{sup}$.  
We {solve  for $\varphi^{sup}$ by using}
the Euler system with the boundary condition at the outer boundary of the domain ( the entering flow). 
For example, in the one-dimensional radial flow model, we use the radial supersonic solution that we have found in Subsection~\ref{sec_exact}.
Hence we have radial velocity $-v=\varphi^{sup}_r=V^{sup}(k_0,r)$, density
$\rho^{sup}(k_0,r)$, pressure $\psup = k_0(\rho^{sup}(k_0,r))^\gamma$ and
the total
inward flux $2 \pi M_0= 2 \pi \rho_0 v_0 R$. Without loss of generality we let
$\varphi^{sup}=0$ at $r=R$.

Let the supersonic flow be converging, that is, the flow moves inward toward the origin, starting from $r=R$ hits a
smooth inner boundary 
denoted by $\Gamma$, which can be considered as the circle $r=a$ in Section~\ref{sec_model}.
If there is no shock,
we assume that $V^{sup}(k_0,r)$ remains supersonic before it hits $\Gamma$.

Let $\Omega$ be the domain in $\mathbb{R}^N$ for an arbitrary $N\ge 1$ enclosed by $\Gamma$ and
$A_R \equiv  B_R \setminus \overline{\Omega}$.
Define the non-empty closed convex set 
\[
{\cal{K}} := \{ w \in W^{1,\infty}(A_R): w \geq \varphi^{sup},\;  |\nabla w |\leq {\hat{q}_0} \}
\]
and  {$c^2 :=  \frac{\gamma-1}{2} (\hat{q}_0^2-|\nabla w|^2)$ as in \eqref{soundspeed}.}
(This is the sound
speed as given by the Bernoulli's law {if $w$ is the velocity potential for the flow; but at this point
$w$ is just any element in ${\cal K}$ and $c^2$
becomes a definition.)}
We introduce the functional $I: {\cal K} \to {\mathbf R}$ such that
\begin{eqnarray} \label{var_functional}
{I(w)} &=& \int_\Omega - (c^2)^{\gamma/(\gamma-1)}  dx \, dy -\int_{A_R\setminus\Omega} \gamma^{\gamma/(\gamma-1)} k_a^{1/(\gamma-1)} \psup dx \, dy \nonumber  \\ 
& & \qquad +   \int_{\Gamma} \gamma^{\gamma/(\gamma-1)} k_a^{1/(\gamma-1)} m(s)\,{w} \, ds,
\end{eqnarray}
where $\Omega \equiv \{ {\bf x} \in A_R: 
{w}({\bf x}) > \varphi^{sup}(\bf{x}) \}$, $m(s)$ is a given smooth
function on $\Gamma$ satisfying $\int_{\Gamma} m(s) ds= 2 \pi M_0$, and
$k_a>k_0$ is a positive constant to be specified later. (In the one dimensional case, $k_a$ is
given by (\ref{k_a}).) 
Define the coincidence set as $A_R \setminus \Omega$.
Let $\varphi$ be a critical point of the functional $I$.
We now derive the jump condition across $\Sigma$, which is the part of
$\partial (A_R \setminus \Omega)$ that lies in the interior of $A_R$,  and study its relation to 
the Rankine-Hugoniot jump conditions. It is easy to see that there is no jump in the  tangential derivatives of $\varphi$ and $\varphi^{sup}$
across $\Sigma$. For the one dimensional case we confirm later that the jump conditions imply
the rest of the Rankine-Hugoniot conditions as well.

{We use  $\partial/\partial n$ to denote the normal derivative in the increasing $r$ direction on the shock $\Sigma$,
and in the outward normal
direction for the remaining two boundaries $\partial B_R$ and $\Gamma$.}
Taking the Fr\'echet derivative, we have
\begin{eqnarray*}
\delta I &=&  \int_\Omega -\frac{\gamma}{\gamma-1}  (c^2)^{1/(\gamma-1)} (-\frac{\gamma-1}{2}) 2 \nabla \varphi \cdot  \nabla \delta\varphi \, dx dy
+\int_{\Sigma} \frac{\delta\varphi}{\partial \phisup/\partial n  -\partial \varphi/\partial n} (-(c^2)^{\gamma/(\gamma-1)}) ds\\
&&+\int_{\Sigma} \frac{\delta\varphi}{\partial \phisup/\partial n -\partial \varphi/\partial n} \gamma^{\gamma/(\gamma-1)} k_a^{1/(\gamma-1)} \psup  \,ds+ \int_{\Gamma} \gamma^{\gamma/(\gamma-1)} k_a^{1/(\gamma-1)} m(s) \delta\varphi  \,ds\;.
\end{eqnarray*}
\begin{figure}
\begin{center}
\includegraphics[height = 2in,width = 2.8in]{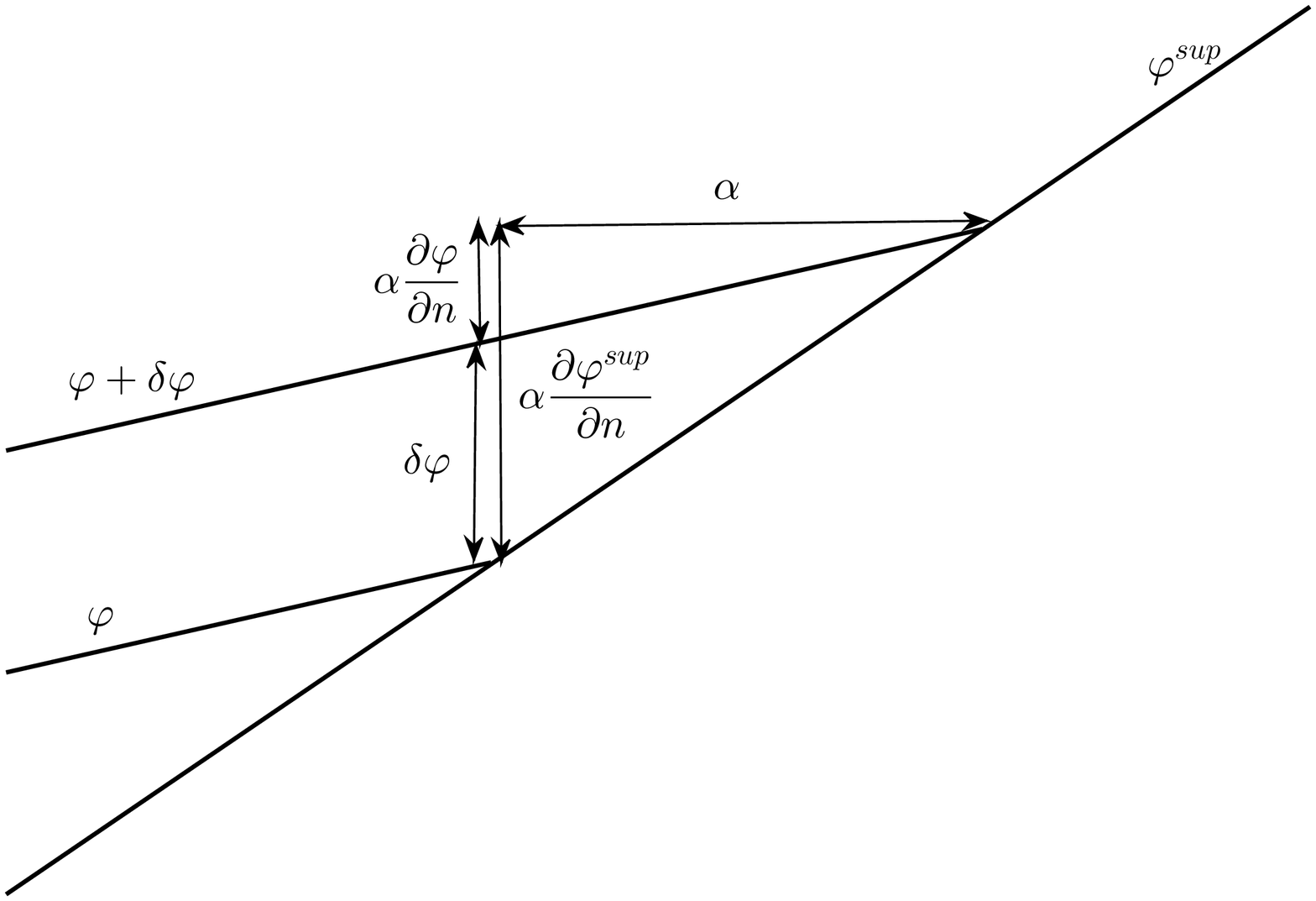}
\end{center}
\caption{The displacement $\alpha$ of the shock position when $\varphi$ changes to $\varphi+\delta\varphi$}
\label{freebc}
\end{figure}
Figure~\ref{freebc} depicts the displacement $\alpha$ of the shock positions {when changing} $\varphi$ to $\varphi+\delta\varphi$.
{To leading order  we obtain}
\begin{eqnarray*}
\alpha &=& \frac{\delta\varphi}{ \partial\varphi^{sup}/\partial n - \partial\varphi/\partial n }.
\end{eqnarray*}

\begin{remark}
There is a well known change of domain formula \cite[Theorem 1.11, p.14]{Henry}
\[
\frac{d}{dt} \int_{\Omega(t)} f(t,x) \,dx dy= \int_{\Omega(t)} \frac{\partial f}{\partial t} \,dx dy+ \int_{\partial \Omega (t)} f(t,x) {\bf V} \cdot {\bf N} ds \;,
\]
where ${\bf V} \cdot {\bf N}$ is the velocity component of the evolving boundary in the normal direction ${\bf N}$. By setting
 ${\bf V} \cdot {\bf N} \, \delta t=\pm \frac{1}{ \partial\varphi^{sup}/\partial n - \partial\varphi/\partial n } \delta \varphi$, 
depending on whether we are dealing with $\Omega$ or $A_R\setminus \Omega$, we can also derive the same formula for $\delta I$.
We also refer Courant and Hilbert \cite[p.260-262]{CoHi} for variable domains for the calculus of variations.
\end{remark}

Assuming a certain smoothness of $\varphi$ and the shock $\Sigma$,
 the first term in the last equation becomes
\begin{eqnarray*}
 \int_\Omega \gamma  (c^2)^{1/(\gamma-1)} \nabla \varphi \cdot \nabla \delta\varphi  \, dx dy
&=&  - \int_\Omega \gamma  \, div ( (c^2)^{1/(\gamma-1)} \nabla \varphi ) \delta\varphi \, dx dy\\
&&+  \int_{\Sigma} \gamma  (c^2)^{1/(\gamma-1)} \frac{\partial \varphi}{\partial n} \delta\varphi \, ds
+ \int_{\Gamma} \gamma  (c^2)^{1/(\gamma-1)} \frac{\partial \varphi}{\partial n}  \delta\varphi  \, ds\;.
\end{eqnarray*}
In $\Omega$, we define $\rho$ by 
\begin{equation} \label{var_rho}
(c^2)^{1/(\gamma-1)} = (k_a \gamma)^{1/(\gamma-1)} \rho
\end{equation}
{which is inspired by (\ref{soundspeed}), and  $p=k_a \rho^{\gamma}$;} hence 
\begin{equation}
(c^2)^{\gamma/(\gamma-1)} = (k_a \gamma)^{\gamma/(\gamma-1)} \rho^\gamma = k_a^{1/(\gamma-1)} \gamma^{\gamma/(\gamma-1)} p \;. \label{var_p} 
\end{equation}
We can now rewrite $\delta I$ as
\begin{eqnarray}
\delta I &=&  - \int_\Omega \gamma  (k_a \gamma)^{1/(\gamma-1)} \,
div (\rho \nabla \varphi) \delta\varphi \, dx dy \nonumber \\
&&+  \int_{\Sigma} \gamma  (k_a \gamma)^{1/(\gamma-1)} \rho  \frac{\partial \varphi}{\partial n} \delta\varphi   \, ds
+ \int_{\Gamma} \gamma   (k_a \gamma)^{1/(\gamma-1)} \rho \frac{\partial \varphi}
{\partial n}  \delta\varphi \, ds  \nonumber \\
&&+\int_{\Sigma} \frac{\delta\varphi}{\partial \phisup/\partial n -\partial \varphi/\partial n} (-k_a^{1/(\gamma-1)} \gamma^{\gamma/(\gamma-1)} p  ) \, ds 
\label{Iprime} \\
&&+\int_{\Sigma} \frac{\delta\varphi}{\partial \phisup/\partial n -\partial \varphi/\partial n} \gamma^{\gamma/(\gamma-1)} k_a^{1/(\gamma-1)} \psup \, ds +  \int_{\Gamma} \gamma^{\gamma/(\gamma-1)} k_a^{1/(\gamma-1)} m(s) \delta\varphi \, ds\;. \nonumber
\end{eqnarray} 
{Suppose the critical point $\varphi$ satisfies $|\nabla \varphi|< \hat{q}_0$ in $\Omega$ so that $\delta \varphi$ can be arbitrary there 
as long as it is small enough,}
we recover \eqref{mass_div}, i.e. $div (\rho \nabla \varphi)=0$ in  $\{{\mathbf{x}} \in A_R: \varphi>\phisup, \;|\nabla \varphi|< {\hat{q}_0} \}$. Moreover,  
\begin{equation} \label{var_bc_a}
\rho \partial \varphi/\partial n= - m(s) \;\mbox{ on} \; \Gamma 
\end{equation}
and
\[
\rho \frac{\partial \varphi}{\partial n}- \frac{p}{\partial \phisup/\partial n -\partial \varphi/\partial n}+\frac{p^{sup}}{\partial \phisup/\partial n -\partial \varphi/\partial n}=0
\quad \mbox{on} \; \Sigma,
\]
which simplifies to
 \begin{equation} \label{jump_cond}
p+ \rho (\frac{\partial \varphi}{\partial n})^2 = p^{sup}+ \rho \frac{\partial \varphi}{\partial n}
\frac{\partial \phisup}{\partial n}   \quad \mbox{on} \; \Sigma \;.
 \end{equation}
This is a linear combination of the remaining two Rankine-Hugoniot conditions
\eqref{sk1} and \eqref{sk2}.

Assume there is only one {continuous} shock {that encloses $\Gamma$,}
by integrating $div (\rho \nabla \varphi)=0$
we have
\begin{equation} \label{flux_conserve}
\int_{\Sigma}\rho \frac{\partial \varphi}{\partial n} 
=-\int_{\Gamma} \rho \frac{\partial \varphi}{\partial n}=
\int_{\Gamma} m(s) ds= 2\pi M_0
=
\int_{\Sigma}\rho^{sup} \frac{\partial \varphi^{sup}}{\partial n} \;.
\end{equation}
Hence the total flux will be conserved.
We have now established the following theorem.
\begin{theorem} \label{main}
Assume the existence and smoothness of a
 critical point $\varphi \in \mathcal{K}$ for the functional $I$ in \eqref{var_functional} and the coincidence set (where $\varphi=\varphi^{sup}$)
 has a smooth boundary $\Sigma$. Then
 the critical point $\varphi$ satisfies 
\begin{eqnarray*}
div (\rho(|\nabla \varphi|^2)\nabla \varphi) &=&0 \quad \mbox{in}\quad  \{{ {\bf x} \in A_R}: \varphi >\varphi^{sup}, \ |\nabla \varphi|< {\hat{q}_0}\} \;,\\
p+ \rho (\frac{\partial \varphi}{\partial n})^2& =& p^{sup}+ \rho \frac{\partial \varphi}{\partial n}
\frac{\partial \phisup}{\partial n}   \quad \mbox{on} \quad \Sigma \;,  \\ 
\rho \partial \varphi/\partial n &=& - m(s) \quad \;\mbox{ on} \quad \Gamma \;, 
\end{eqnarray*}
where $m$ is a given function that satisfies (\ref{flux_conserve}).
\end{theorem}

\section{Characterization of the critical point} \label{sec_analysis1D}
\setcounter{equation}{0}

For the model problem, let $\Gamma=\partial B_a$, and impose the same boundary conditions (\ref{bcR}) and (\ref{bca}) for the variational problem with a subsoinc $v_a$ 
and a supersonic $v_0$. Thus $M_0=\rho_0 v_0 R$ is known. Set
  $m=M_0/a$ so that (\ref{flux_conserve}) is satisfied. With $\rho(a)=m/v_a$, $k_a$ can be calculated from (\ref{var_rho}):
 \[
k_a=
\frac{\gamma -1}{2\gamma} \left(\frac{a v_a}{M_0}\right)^{\gamma-1}(\qzero - v_a^2).
\]
We note that for general multi-dimensional configurations, 
the constant $k_a$ may not be determined uniquely, and may be with the same constant $k_0=k_a$ in the entire region when the shock is weak. 

 Moreover we assume that  the prescribed $v_a$ will make
 $k_a$ satisfy \eqref{k_condition}. 

We need to ensure that $V^{sup}(k_0,r)$ remains supersonic
for $a \leq r \leq R$. 
As $r$ decreases from $R$,
$\rho |v|$ has to increase to maintain a constant flux $M_0$. In the supersonic
regime, $|v|$ has to decrease for $\rho |v|$ to increase. Once $|v|$ reaches 
the sonic speed $c_*$ before reaching $r=a$, the local flux $\rho |v|$ cannot increase
further and there is no supersonic flow afterward. Hence 
in order that $\phisup$ to be defined on the entire region $a \leq r \leq R$,
we impose
\[
\rho_* c_* a >M_0  \;,
\]
where $\rho_*$ is the density at the sonic speed.
Since $\rho_*=(\frac{c_*^2}{\gamma k_0})^{1/(\gamma-1)}$, the above
relation is equivalent to
\begin{equation} \label{condition1}
k_0 < \frac{1}{\gamma} \left( \frac{a}{M_0} \right)^{\gamma-1} c_*^{\gamma+1} \;,
\end{equation}
which is already satisfied by imposing condition \eqref{k_condition}.

If we use $dx \, dy=r dr d\theta$
in the derivation in Section~\ref{sec_variation}, or simply convert the equations to polar coordinates at the end,
the governing transonic flow equation for radially symmetric solutions {in Theorem~\ref{main}}   becomes
\begin{eqnarray}\label{mass1}
(\rho \varphi_r r)_r&=&0
\end{eqnarray}
in $\{ {\bf x} \in A_R: \varphi>\phisup, \; |\varphi_r|< {\hat{q}_0} \}$.
Since $\varphi_r(a)=v_a<c_*$,  
\eqref{mass1} holds in a neighborhood
of $r=a$ with $\rho \varphi_r r= M_0$. Thus $\varphi_r>0$ and
as $r$ increases, $\rho \varphi_r=M_0/r$
decreases. It is easy to show that $\rho \varphi_r$ is an increasing function of 
the speed $\varphi_r$ so long as $\varphi_r<c_*$, and thus $\varphi_r$ decreases with larger $r$ and will never reach $c_*$.
{Hence the equation $\rho \varphi_r=M_0/r$ can uniquely be solved to get  a subsonic $\varphi_r$ in term of $M_0$, $k_a$ and $r$.}

We now locate the shock by a variational method.
 Let  $a \leq \eta \leq R$. Define
\begin{equation}  \label{K1}
{\cal{K}}_1 \equiv 
\left\{ \begin{array}{ll}
 \varphi \in W^{1,\infty}(A_R): 
 & \varphi = \varphi^{sup} \; \mbox{for} \; \eta \leq r \leq R, \\ 
&  \rho \varphi_r r=M_0 \; \mbox{for} \; a \leq r \leq \eta, \\ 
& \mbox{and}\; \varphi(\eta)=\varphi^{sup}(\eta) 
\end{array}
\right\}
\end{equation}
With a given $\eta$, the function $\varphi \in {{\cal K}}_1$ is uniquely determined.
We now consider  $\varphi \in {{\cal K}}_1$  as a function of $\eta$ and
define $J: (a,R) \to \bf{R}$ such that $J(\eta)=I(\varphi(\eta))$ with $\varphi \in
{{\cal K}}_1$.

{Since $\varphi \in C^2[a,\eta]$ satisfies the transonic flow equation and the boundary condition at $r=a$, 
(\ref{Iprime})
can be simplified to}
\[
\delta J= \gamma^{\gamma/(\gamma-1)} k_a^{1/(\gamma-1)}
\frac{\delta \varphi}{\varphi^{sup}_r-\varphi_r} 2 \pi \eta \{
\rho \varphi_r (\varphi^{sup}_r-\varphi_r) - p + p^{sup} \} \;.
\] 
Recall that  $\delta r=\delta \varphi/(\varphi^{sup}_r-\varphi_r)$.
Thus we can evaluate $J'$ explicitly to
\begin{equation} \label{J_prime}
J'(\eta)=2 \pi  \gamma^{\gamma/(\gamma-1)} k_a^{1/(\gamma-1)}
\eta \{
p^{sup}+ \rho \varphi_r \varphi^{sup}_r - p - \rho \varphi_r^2  \}|_{r=\eta} \;.
\end{equation}
If an interior critical point of $J$ exists at $\eta=r_s$,  then
 \begin{equation} \label{jump_cond1D}
 \rho \varphi_r^2 + p = \rho \varphi_r\phisup_r +\psup \;\quad \mbox{at} \; r=r_s.
 \end{equation}

As $m$ is prescribed to ensure (\ref{flux_conserve}) is satisfied, thus
$\rho\varphi_r = \rho^{sup} \phisup_r$ at $r=r_s$. Combining this equation with equation \eqref{jump_cond1D}, we have
 \begin{eqnarray*}
 \rho \varphi_r^2 + p = \rho^{sup}( \phisup_r)^2 +\psup \;.
 \end{eqnarray*}
Hence all the Rankine-Hugoniot conditions are satisfied at the shock.

In fact by using  $\varphi \in {{\cal K}}_1$ and
the same calculations 
in obtaining \eqref{flux_conserve},
we always have $\rho \varphi_r= \rho^{sup}
\varphi^{sup}_r$ at any $\eta$, which may  not be a critical point of $J$. Hence
with the definition \eqref{def_H}, we can rewrite \eqref{J_prime} as
\begin{equation} \label{J_prime1}
J'(\eta)=2 \pi  \gamma^{\gamma/(\gamma-1)} k_a^{1/(\gamma-1)}
 (H^{sup}(k_0,\eta) - H^{sub}(k_a,\eta))
 \;.
\end{equation}
We already known that $dH^{sub}/d \eta>dH^{sup}/d\eta$ in Section~\ref{sec_exact},
and thus an interior critical point of $J$, if it exists, is unique.

To show there exists an interior critical point of $J$, it suffices to assume 
that $J'(a^+)>0$ and $J'(R^-)<0$. In other words
\begin{eqnarray*}
J'(R^-) &=& 2 \pi k_a^{1/(\gamma-1)} \gamma^{\gamma/(\gamma-1)} [H^{sup}(R^-) -H^{sub}(R^-)]=2 \pi k_a^{1/(\gamma-1)} \gamma^{\gamma/(\gamma-1)} [ H_C-H_D]<0,\\
J'(a^+) &=&2 \pi k_a^{1/(\gamma-1)} \gamma^{\gamma/(\gamma-1)} [ H^{sup}(a^+) -H^{sub}(a^+)]= 2 \pi k_a^{1/(\gamma-1)} \gamma^{\gamma/(\gamma-1)} [H_B-H_A]>0.
\end{eqnarray*}
They are the same as the necessary and sufficient conditions (\ref{cond1})-(\ref{cond2}).
 This unique interior critical point is a local maximum for $J$.
 
 Now suppose $J'(a^+)<0$ and $J'(R^-)>0$. This leads to $H_C>H_D$ and $H_A>H_B$. As $H^{sup}$ and
 $H^{sub}$ are increasing functions, we have $H_C>H_D>H_A>H_B$. It contradicts $dH^{sub}/d \eta>dH^{sup}/d\eta$.
 This case can never happen.

If 
$\delta \varphi\equiv w=0$ at the shock $\Sigma$, 
the term involving $\int_{A_R \setminus \Omega}$  in (\ref{var_functional})
does not change and the subsonic domain is fixed. 
{It is known that the solution of the transonic equation in the subsonic case can be obtained as a minimizer of $I$, see \cite{DongOu}
which has essentially the same functional.}
We verify {that the subsonic solution is the minimizer of $I$} for our model.

Recall that $dx dy$ becomes $2 \pi r dr$ and 
$ds$ on $\Gamma$ becomes $2 \pi a$, hence
\[
I'(\varphi)w=2 \pi \int_a^\eta \gamma  (c^2)^{1/(\gamma-1)} \varphi_r  w_r r \, dr
 +2 \pi \gamma^{\gamma/(\gamma-1)} k_a^{1/(\gamma-1)} M_0 w(a) \;,
\]
with critical point $\varphi$ satisfying $(\rho \varphi_r r)_r=0$
with boundary conditions $\rho \varphi_r a =M_0$ at $r=a$, and $\varphi=\varphi^{sup}$
at $r=\eta$.
Taking another derivative,
\begin{eqnarray*}
I''(\varphi) (w,w)&=&2\pi \gamma \int_a^\eta \{ - (c^2)^{(2-\gamma)/(\gamma-1)}
(\varphi_r   w_r )^2 + (c^2)^{1/(\gamma-1)} w_r^2 \} r\\
&{=} & 2\pi \frac{\gamma (\gamma+1)}{2} \int_a^\eta (c^2)^{(2-\gamma)/(\gamma-1)}w_r^2
 ( c_*^2 - \varphi_r^2 ) r,
\end{eqnarray*}
which shows that  $I$ is strictly convex for subsonic flow. Hence when we varies
$\varphi$ with $\delta \varphi=0$ at $\Gamma$ (so that the subsonic
domain $[a,\eta]$ is fixed), $I$ attains its minimum when
$\varphi$ satisfies the transonic equation for fixed $\eta$.
This accounts for our choice of the set ${\cal K}_1$. Thus we expect
the critical point that we are looking for may be characterized as a saddle point by
\[
\displaystyle \sup_{a<\eta<R} \{
\min_{\begin{array}{c} \varphi \in {\cal K}, \\ shock \;at \; r=\eta
\end{array}} I(\varphi) \} = \sup_{a<\eta<R} J(\eta) \;.
\]
{Because we can solve the 1D variational problem exactly,} 
we have the following saddle node theorem. 
 \begin{theorem}
 There exists a unique transonic solution $\varphi\in{{\cal K}}_1\subset W^{1,\infty}([a,R])$ to
 \begin{eqnarray*}
 \sup_{a<\eta<R} \{
\min_{\begin{array}{c} \varphi \in {\cal K}_1, \\ shock \;at \; r=\eta
\end{array}} I(\varphi) \}   = \sup_{a<\eta<R} J(\eta) \;,
 \end{eqnarray*} 
 if and only if
 \begin{eqnarray*}
 J'(R^-) <0, \  & and &\  J'(a^+)>0.
 \end{eqnarray*}
 \end{theorem}

\begin{remark}
 The functional $I$ is well defined when we restrict ourselves to
 $W^{1,\infty}$ functions. To prove existence of a critical point, a cut-off of our functional
may be introduced (such as the one in \cite{DongOu}) so that a different function space can be employed as its domain.

The natural question is whether the procedures developed in this paper can be extended to more general setting of the multi-dimensional case. We leave the following open question for future work.
{\em Given any curve ${\cal C}$ in the interior of the flow domain
$A_R$., we find a minimizer $\varphi_{\cal C}$ of the "subsonic" problem. Let $J({\cal C})=I(\varphi_{\cal C})$. Now identify a 
${\cal C}=\Sigma$ such that $J$ attains either its maximum or a saddle point.}
\end{remark}

\section{Numerical results for {variational formulation} of the model problem} \label{sec_num}
This section comprises the numerical results that validate our variational formulation and the model problem.
In all computations, we use $\gamma=3$, $\cstar=1$, $M_0=1$, $R=6$, and $a=5$. The figures are generated by MATLAB.
As in the last section we let 
\[
J(\eta)=\min_{\begin{array}{c} \varphi \in {\cal K}_1, \\ shock \;at \; r=\eta
\end{array}} I(\varphi) \;.
\]
For clarity, we use $\eta_s$ instead of $\eta$ to denote the arbitrarily assigned shock position; the computed
value of the true position of the shock is denoted by $r_s$.

Since $\gamma=3$, {from the Bernoulli's law and the mass conservation} 
\begin{eqnarray*}
\frac{1}{2} V^2 + \frac{k \gamma}{\gamma -1} \bigg( \frac{M_0}{rV}\bigg)^{\gamma-1} &=& \frac{1}{2} {\hat{q}_0^2} \;,
\end{eqnarray*} 
 we write  $V$ (which we only consider positive) explicitly to 
\begin{eqnarray}\label{Veq}
V_{\pm}(k, r) &:=&V= \frac{1}{2} \sqrt{ {\hat{q}_0^2} \pm \sqrt{ \hat{q}_0^4 - 4 \frac{2 k \gamma}{\gamma -1} \bigg( \frac{M_0}{r}\bigg)^{\gamma-1}} }.
\end{eqnarray}
{Since $c_*=\sqrt{\frac{\gamma-1}{\gamma+1}}  \hat{q}_0=\frac{1}{\sqrt{2}} \hat{q}_0$, it is clear that $V_+=V^{sup}$ is supersonic while
$V_-=V^{sub}$ is subsonic. }

\begin{figure}
\centering
\includegraphics[height = 1.8in,width = 2.3in]{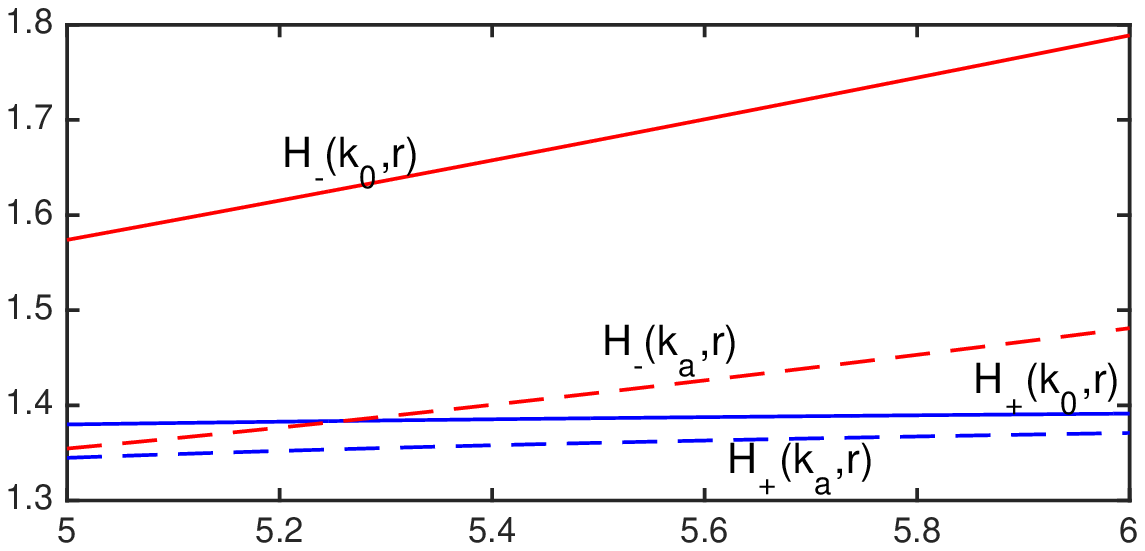}
\includegraphics[height = 1.8in,width = 2.3in]{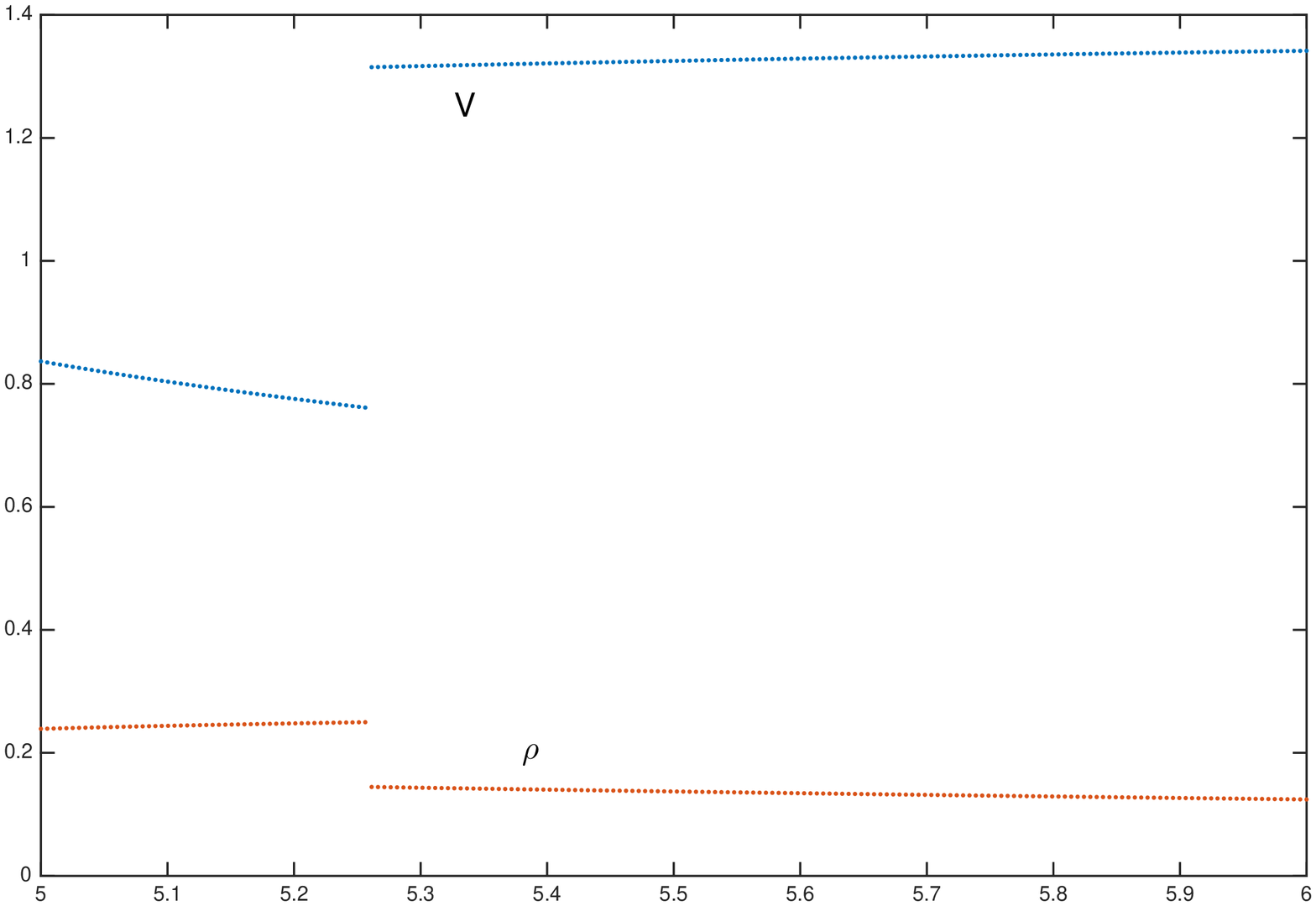}
\vspace*{-4mm}
\caption{Left: $H_{\pm}(k,r)= \frac{M_0 (\gamma+1)}{2\gamma} (\frac{c_*^2}{V_{\pm}(k,r)} + V_{\pm}(k,r))$  when $k=k_a$ (dashed lines) and $k=k_0$ (solid lines). 
Right: The velocity $V$ and the density $\rho$ plots showing the jump at $r_s$.}
\label{HSolns}
\end{figure}

\begin{figure}
\centering
\includegraphics[height = 1.8in,width = 2.3in]{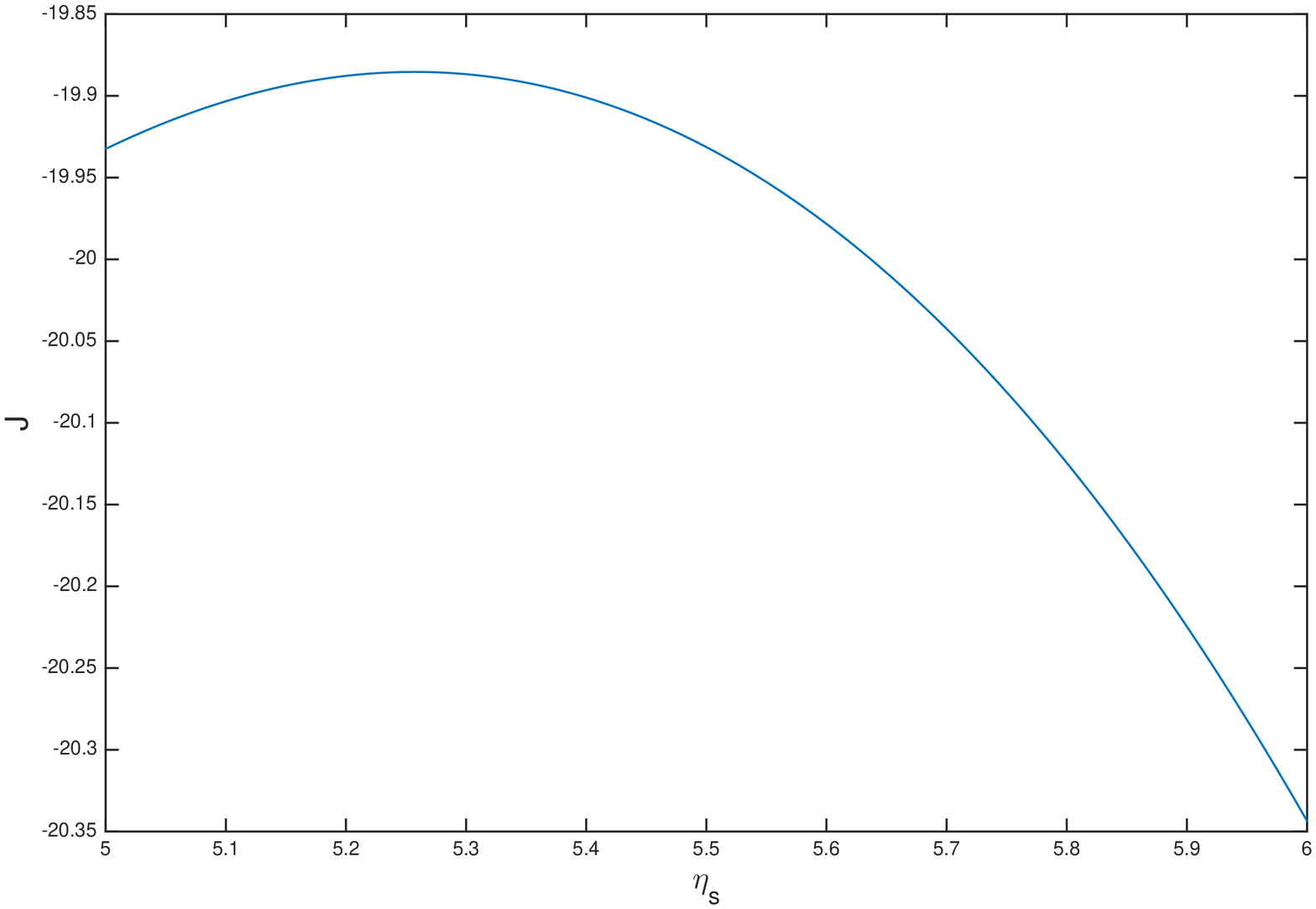}
\includegraphics[height = 1.8in,width = 2.3in]{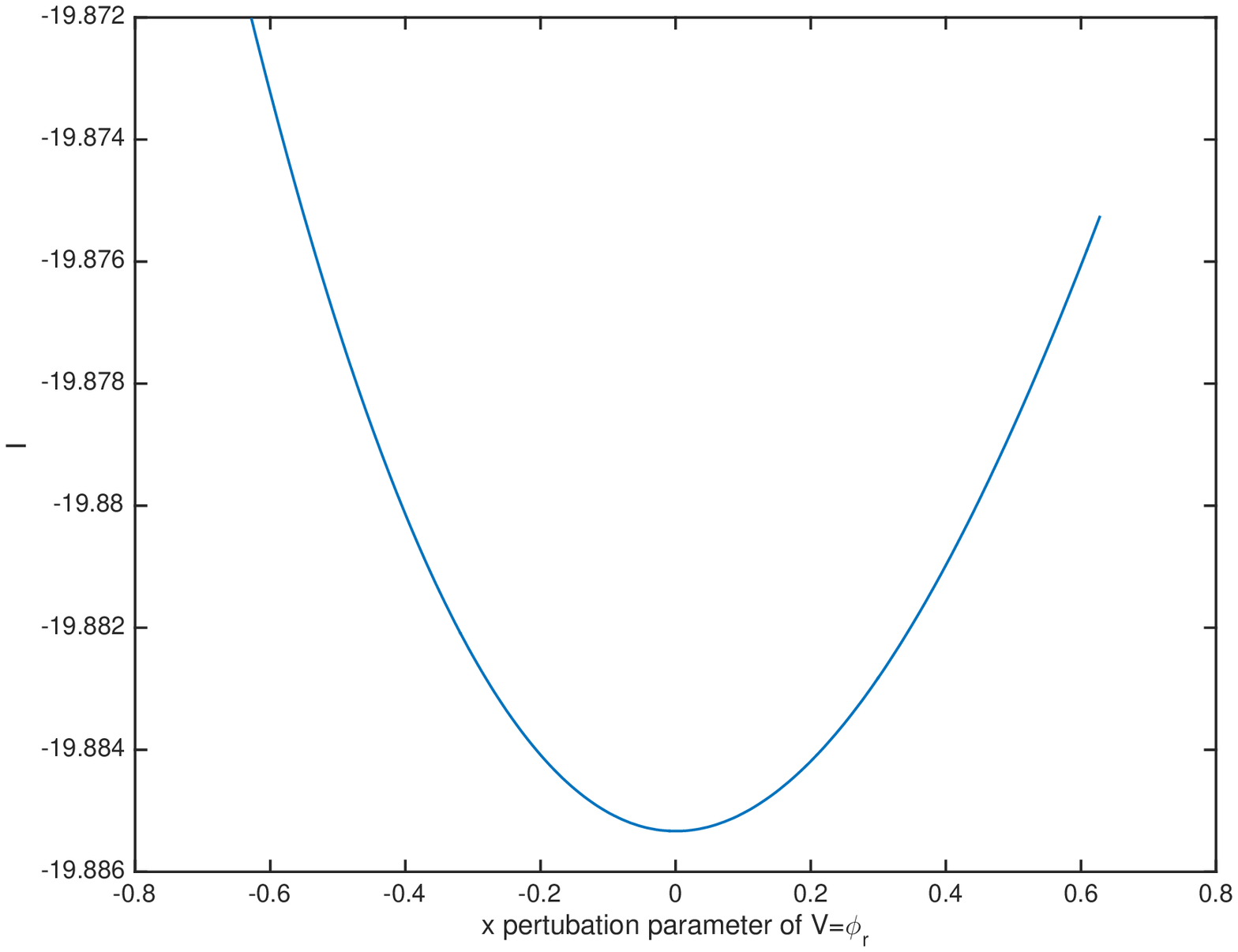}
\vspace*{-4mm}
\caption{Left: {$J/2\pi$} graph showing the maximum at $\eta_s=r_s$. Right: $I(x)/2\pi$ graph showing the minimum at ${x=0}$, i.e. at $V^{sub}$ when $\eta_s=r_s$.}
\label{JIvariationals}
\end{figure}

For Figures~\ref{HSolns},~\ref{JIvariationals}, we have $v_a= \sqrt{0.7}\approx 0.8367 > \nu =\sqrt{0.5}$ and $v_0=\sqrt{1.8}\approx 1.3416$, {where $\nu$ is defined in \eqref{vacond}.} The corresponding $k$ values are $k_a\approx 7.58\overline{3}$ and $k_0\approx 4.32$.
{Our construction ensures that $V_+(k_0,R)=v_0$ and $V_-(k_a,a)=v_a$.}

{The left figure in Figure~\ref{HSolns} depict graphs of $H_{\pm}(k,r)= \frac{M_0 (\gamma+1)}{2\gamma} (\frac{c_*^2}{V_{\pm}(k,r)} + V_{\pm}(k,r))$;
$H_+=H^{sup}$ and $H_-=H^{sub}$ are computed by substituting $V$ by $V_+(k_0,r)$ and $V_-(k_a,r)$, respectively.
As expected there exists a unique intersection point for the graphs of $H^{sup}$ and $H^{sub}$
 at  $\eta_s=r_s\approx 5.260220746$. The resulting velocity $V=V^{sub/sup}$ and the  corresponding density $\rho=\rho^{sub/sup}$ are depicted in the right figure in Figure~\ref{HSolns}. Such a solution satisfies the physical entropy condition,
namely density increases across the shock in the direction of flow.}

{Many studies make an assumption that the entropy is with the same constant, for example $k=k_a=k_0$ in the entire region. }
This may lead to an erroneous conclusion when the shock is strong. 
It can be easily shown that both $H^{sup}(\cdot,r)$ and $H^{sub}(\cdot,r)$ are decreasing functions of $k$.
That explains why  the graph of $H^{sub}$ when $V$ is replaced by $V_-(k_0,r)$,
represented by {the solid line labeled as $H_-(k_0,r)$} in the left figure of Figure~\ref{HSolns},
 is located well above the correct $H^{sub}$.
This curve  does not intersect $H^{sup}$ and give the incorrect conclusion that no shock develops.

Two figures in Figure~\ref{JIvariationals} confirm that the critical point of $I$ is a saddle node.
The left figure in Figure ~\ref{JIvariationals} is the graph of {$J/2\pi$}  with respect to $\eta_s$, where
\begin{eqnarray*}
{\frac{1}{2\pi}} J(\eta_s) &=&  \int_a^{\eta_s} -({c^2 \Big|_{V=V_-(k_a,r)}})^{\gamma/(\gamma-1)} r dr 
+ \int^R_{\eta_s} -\gamma^{\gamma/(\gamma-1)} k_a^{1/(\gamma-1)} {p \Big|_{V=V_+(k_0,r)}}  r dr \\
&&+ \gamma^{\gamma/(\gamma-1)} k_a^{1/(\gamma-1)} M_0 \bigg( \int_a^{\eta_s} -V_-(k_a,r) dr +  \int^R_{\eta_s} -V_+(k_0,r) dr \bigg),
\end{eqnarray*}
which has a unique maximum point at $\eta_s=r_s$. {Here we have employed $M_0=m a$ and $\varphi=0$ at $r=R$.}

The right figure in Figure~\ref{JIvariationals} is the graph of {$I/2\pi$} 
with respect to a perturbation parameter $x$. {Precisely
let $\Delta V= x (r-a)$ be the perturbation in speed in the subsonic region while keeping the assigned shock location $\eta_s$ fixed. 
The given speed at  the boundary $r=a$ is unchanged. One can also assume we have the same $\varphi$ at $\eta_s$, as 
the functional depends on the derivative of $\varphi$. Regarding $I$ as a function of $x$, we have
\begin{eqnarray*}
{\frac{1}{2\pi}} I(x) &=&  \int_a^{r_s} -({c^2 \Big|_{V=V_-(k_a,r)+\Delta V}})^{\gamma/(\gamma-1)} r dr 
+ \int^R_{r_s} - \gamma^{\gamma/(\gamma-1)} k_a^{1/(\gamma-1)} {p \Big|_{V=V_+(k_0,r)}} r dr \\
&&+ \gamma^{\gamma/(\gamma-1)} k_a^{1/(\gamma-1)} M_0 \bigg( \int_a^{r_s} -(V_-(k_a,r)+\Delta V) dr +  \int^R_{r_s} -V_+(k_0,r) dr \bigg),
\end{eqnarray*}
with {$x \in [-0.62,0.62]$} so that $V^2 < c_*^2$ (note that the maximum value of $V=0.83\bar{6}$ in this example).}
The graph shows that 
$I$ attains its minimum at
$V^{sub}$ which {corresponds to} $x=0$. 

\begin{figure}
\centering
\includegraphics[height = 1.8in,width = 2.3in]{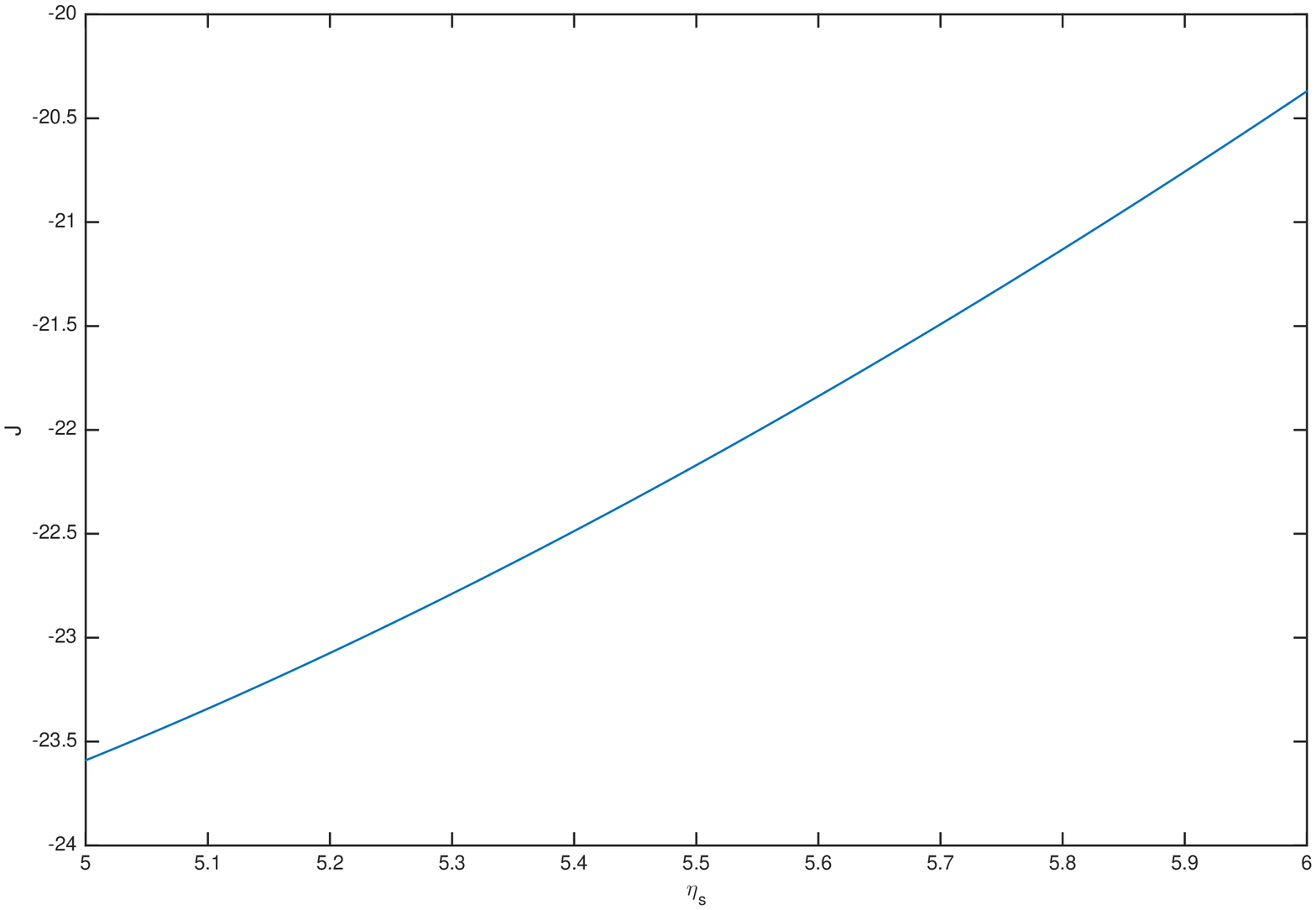}
\includegraphics[height = 1.8in,width = 2.3in]{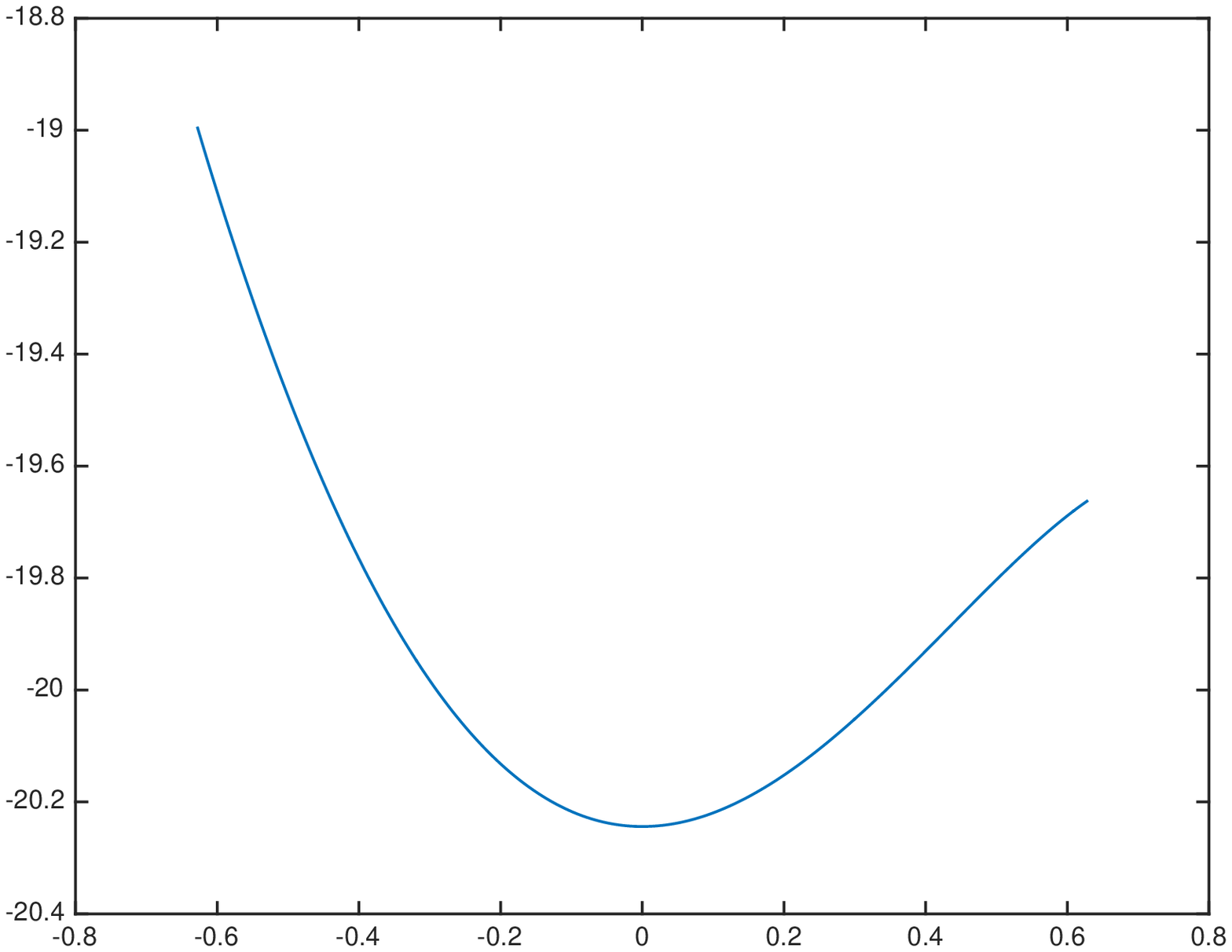}
\vspace*{-4mm}
\caption{Subsonic flow -- Left: {$J/2\pi$} graph showing the maximum at $\eta_s=R$. Right: $I(x)/2\pi$ graph showing the minimum at $x=0$ (that is at $V^{sub}$).}
\label{noshocks}
\end{figure}

\begin{figure}
\centering
\includegraphics[height = 1.8in,width = 2.3in]{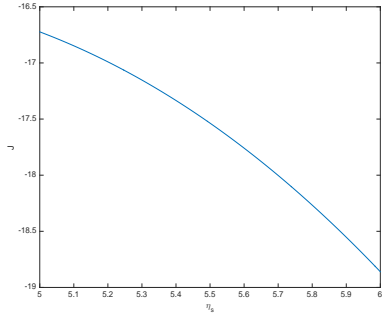}
\vspace*{-4mm}
\caption{Supersonic flow -- {$J/2\pi$} graph showing the maximum at $a$}
\label{cegnoshk}
\end{figure}

{In Figures~\ref{noshocks}, ~\ref{cegnoshk}, we show numerics for the critical point without a transonic shock; there is no transonic flow.}

Two figures in Figures~\ref{noshocks} are for the case when the flow is subsonic in the entire region. For computation, we have used $v_a=\sqrt{0.7}$, and $k_a=k_0\approx 7.58\overline{3}$.
 To capture the subsonic flow in the entire region, the data on $r=R$ is ignored. {The lack of a shock wave allows us to use $k_a=k_0$.}
The left figure is the graph of {$J/2\pi$} showing its maximum at $\eta_s=R$, and the right figure is the graph of $I(x)$ showing its minimum at $V=V^{sub}$ (when it is unperturbed). 

{Figure~\ref{cegnoshk} depicts the graph of  {$J/2\pi$} when the flow is supersonic everywhere.}
This results was discussed earlier in Subsection~\ref{sec_exact}.
In this case we have $v_a= \sqrt{0.4} < \nu =\sqrt{0.5}$ and $v_0=\sqrt{1.8}$. {The maximum of $J$ is attained at $\eta_s=a$ {instead}.}

When a shock is not formed in the interior of the flow domain,
 the {maximum of $J$} is attained either  at $r=a$ or $r=R$ depending on the boundary data. 
 {The corresponding solution $\varphi \in {\cal K}_1$
 does not satisfies the boundary conditions at {either} $r=R$ {or} $r=a$,}
  since it remains either supersonic ({when J attains its maximum} at $r=a$ as in Figure~\ref{cegnoshk}) or subsonic  ({maximum} at $r=R$ as in Figure~\ref{noshocks}) in the entire region.

\begin{figure}
\centering
\includegraphics[height = 1.8in,width = 2.3in]{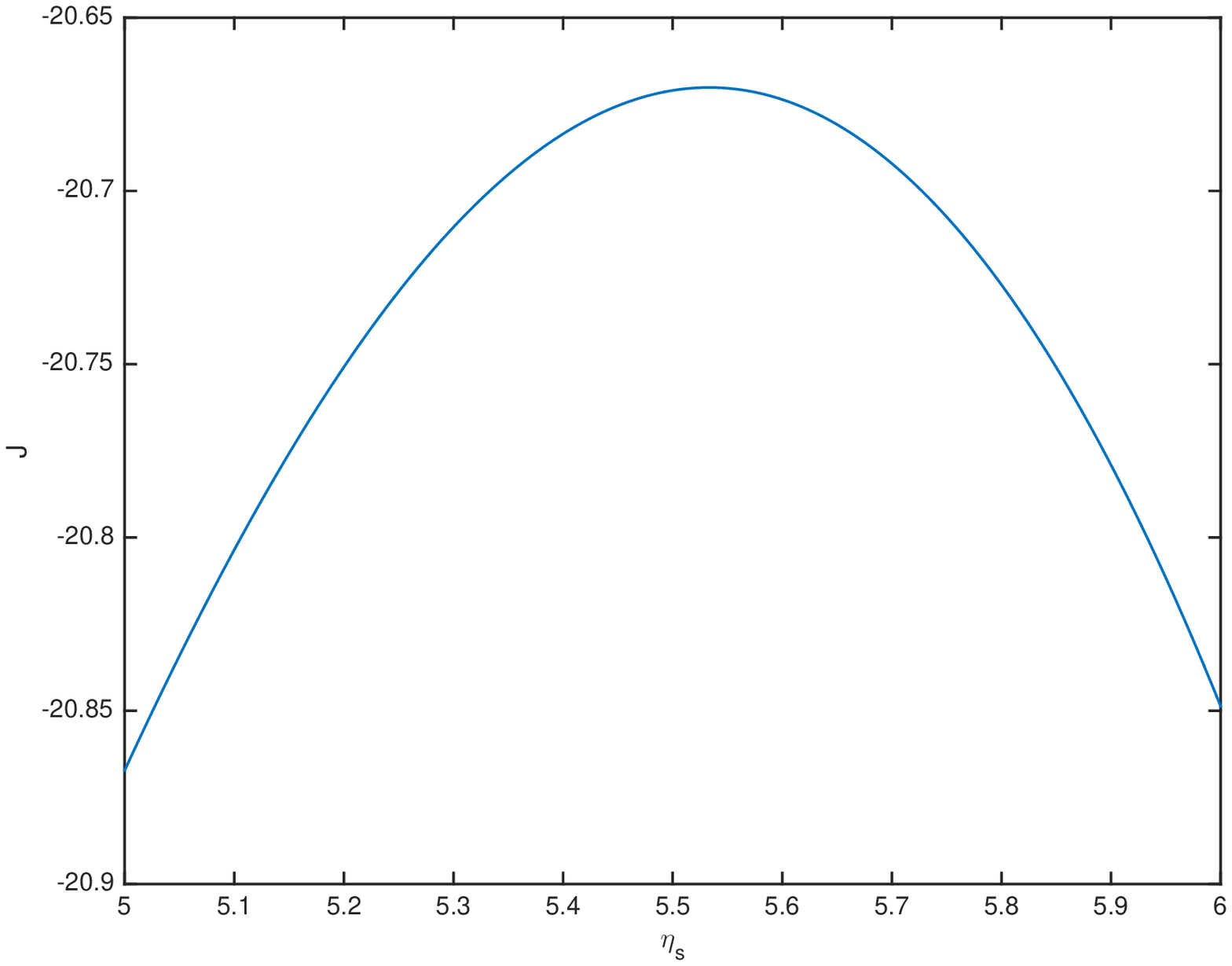}
\includegraphics[height = 1.8in,width = 2.3in]{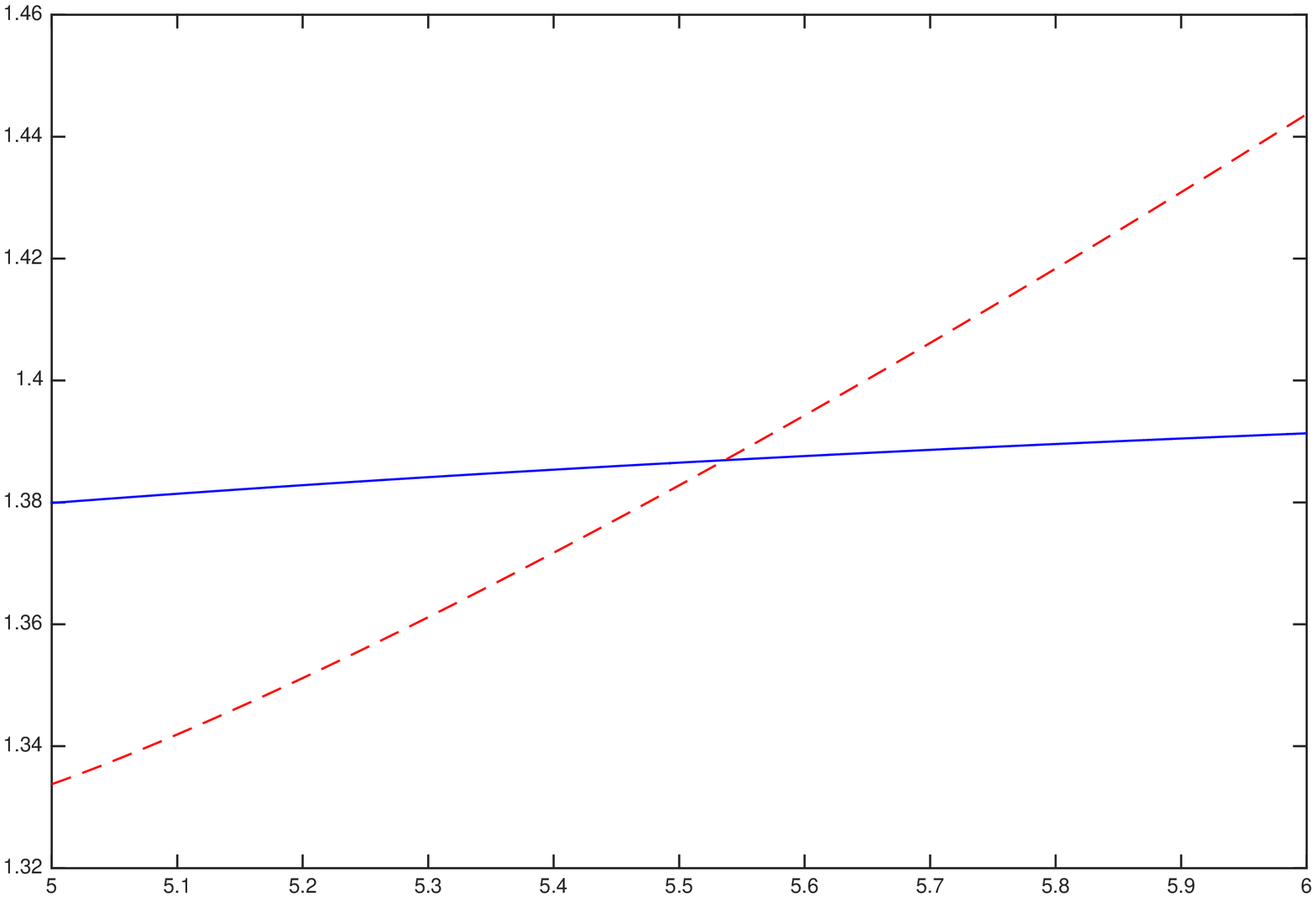}
\vspace*{-4mm}
\caption{Larger $V_a$ --  Left: {$J/2\pi$} graph. 
Right:  $H$ graphs for $k_a$ (dashed curve) and $k_0$ (solid curve).}
\label{eg2}
\end{figure}

Figure ~\ref{eg2} is the case for a larger value of $V_a$ where the values $v_a= \sqrt{0.95}\approx 0.9747> \nu =\sqrt{0.5}$ and $v_0=\sqrt{1.8}$ are used for computations. The corresponding $k$ values are $k_a\approx 8.3125$ and $k_0\approx 4.32$. 
The maximum value of $J$ is attained at $r_s\approx 5.5363$ (see {the right figure of } Figure~\ref{eg2} where $H$ curves are intersecting)
whereas the shock position $r_s=5.260220746$ obtained in the example in Figure~\ref{Hcurves}.

Finally we compare the necessary and sufficient conditions for the transonic shock and the first variation by evaluating them explicitly. 
More precisely, we use the same data as in Figure~\ref{Hcurves} and obtain
\begin{eqnarray*}
{\frac{1}{2\pi} \frac{\partial J}{\partial \eta}}\mid_{\eta=R}&=&-1.28456894\\
k_a^{1/(\gamma-1)} \gamma^{\gamma/(\gamma-1)} [ H_C-H_D]&=&-1.284568972\\
\frac{1}{2 \pi} \frac{{\partial J}}{\partial \eta}\mid_{\eta=a}&=&0.362180988\\
k_a^{1/(\gamma-1)} \gamma^{\gamma/(\gamma-1)} [H_B-H_A] &=&0.3621809759.
\end{eqnarray*}

We conclude the paper with the following remark.
Our variational inequality formulation holds for general multidimensional setting.
The boundary of its coincidence set forms a shock, across which two out of three Rankine-Hugoniot jump conditions are satisfied.
A good estimate of the 
(average) change in the entropy constant in the configuration would be useful especially when a strong shock develops.
However even if we set $k_a=k_0$ the formulation still stands provided the existence of a critical point.
We leave the further study on general multidimensional transonic problems in variational inequality formulations for future work.

\end{document}